\documentclass[twoside,leqno]{article}
\usepackage{amsmath}
\usepackage{amsthm}
\usepackage{amssymb}
\usepackage{amscd}

\usepackage{latexsym}
\usepackage{amsfonts}

\input xy
\xyoption{all} \tolerance=500

=1in \textwidth 15.6cm \topmargin 15mm\textheight23cm \advance\topmargin by -\headheight\advance\topmargin by
-\headsep

\numberwithin{equation}{section} \allowdisplaybreaks

\theoremstyle{definition}
\newtheorem{definition}{Definition}[section]

\begin{document}
\font\black=cmbx10 \font\sblack=cmbx7 \font\ssblack=cmbx5 \font\blackital=cmmib10  \skewchar\blackital='177
\font\sblackital=cmmib7 \skewchar\sblackital='177 \font\ssblackital=cmmib5 \skewchar\ssblackital='177
\font\sanss=cmss10 \font\ssanss=cmss8 
\font\sssanss=cmss8 scaled 600 \font\blackboard=msbm10 \font\sblackboard=msbm7 \font\ssblackboard=msbm5
\font\caligr=eusm10 \font\scaligr=eusm7 \font\sscaligr=eusm5 \font\blackcal=eusb10 \font\fraktur=eufm10
\font\sfraktur=eufm7 \font\ssfraktur=eufm5 \font\blackfrak=eufb10

\font\bsymb=cmsy10 scaled\magstep2
\def\all#1{\setbox0=\hbox{\lower1.5pt\hbox{\bsymb
       \char"38}}\setbox1=\hbox{$_{#1}$} \box0\lower2pt\box1\;}
\def\exi#1{\setbox0=\hbox{\lower1.5pt\hbox{\bsymb \char"39}}
       \setbox1=\hbox{$_{#1}$} \box0\lower2pt\box1\;}

\def\mi#1{{\fam1\relax#1}}
\def\tx#1{{\fam0\relax#1}}

\newfam\bifam
\textfont\bifam=\blackital \scriptfont\bifam=\sblackital \scriptscriptfont\bifam=\ssblackital
\def\bi#1{{\fam\bifam\relax#1}}

\newfam\blfam
\textfont\blfam=\black \scriptfont\blfam=\sblack \scriptscriptfont\blfam=\ssblack
\def\rbl#1{{\fam\blfam\relax#1}}

\newfam\bbfam
\textfont\bbfam=\blackboard \scriptfont\bbfam=\sblackboard \scriptscriptfont\bbfam=\ssblackboard
\def\bb#1{{\fam\bbfam\relax#1}}

\newfam\ssfam
\textfont\ssfam=\sanss \scriptfont\ssfam=\ssanss \scriptscriptfont\ssfam=\sssanss
\def\sss#1{{\fam\ssfam\relax#1}}

\newfam\clfam
\textfont\clfam=\caligr \scriptfont\clfam=\scaligr \scriptscriptfont\clfam=\sscaligr
\def\cl#1{{\fam\clfam\relax#1}}

\newfam\frfam
\textfont\frfam=\fraktur \scriptfont\frfam=\sfraktur \scriptscriptfont\frfam=\ssfraktur
\def\fr#1{{\fam\frfam\relax#1}}

\def\cb#1{\hbox{$\fam\gpfam\relax#1\textfont\gpfam=\blackcal$}}

\def\hpb#1{\setbox0=\hbox{${#1}$}
    \copy0 \kern-\wd0 \kern.2pt \box0}
\def\vpb#1{\setbox0=\hbox{${#1}$}
    \copy0 \kern-\wd0 \raise.08pt \box0}

\def\pmb#1{\setbox0\hbox{${#1}$} \copy0 \kern-\wd0 \kern.2pt \box0}
\def\pmbb#1{\setbox0\hbox{${#1}$} \copy0 \kern-\wd0
      \kern.2pt \copy0 \kern-\wd0 \kern.2pt \box0}
\def\pmbbb#1{\setbox0\hbox{${#1}$} \copy0 \kern-\wd0
      \kern.2pt \copy0 \kern-\wd0 \kern.2pt
    \copy0 \kern-\wd0 \kern.2pt \box0}
\def\pmxb#1{\setbox0\hbox{${#1}$} \copy0 \kern-\wd0
      \kern.2pt \copy0 \kern-\wd0 \kern.2pt
      \copy0 \kern-\wd0 \kern.2pt \copy0 \kern-\wd0 \kern.2pt \box0}
\def\pmxbb#1{\setbox0\hbox{${#1}$} \copy0 \kern-\wd0 \kern.2pt
      \copy0 \kern-\wd0 \kern.2pt
      \copy0 \kern-\wd0 \kern.2pt \copy0 \kern-\wd0 \kern.2pt
      \copy0 \kern-\wd0 \kern.2pt \box0}

\def\cdotss{\mathinner{\cdotp\cdotp\cdotp\cdotp\cdotp\cdotp\cdotp
        \cdotp\cdotp\cdotp\cdotp\cdotp\cdotp\cdotp\cdotp\cdotp\cdotp
        \cdotp\cdotp\cdotp\cdotp\cdotp\cdotp\cdotp\cdotp\cdotp\cdotp
        \cdotp\cdotp\cdotp\cdotp\cdotp\cdotp\cdotp\cdotp\cdotp\cdotp}}

\font\frak=eufm10 scaled\magstep1 \font\fak=eufm10 scaled\magstep2 \font\fk=eufm10 scaled\magstep3
\font\scriptfrak=eufm10 \font\tenfrak=eufm10


\mathchardef\za="710B  
\mathchardef\zb="710C  
\mathchardef\zg="710D  
\mathchardef\zd="710E  
\mathchardef\zve="710F 
\mathchardef\zz="7110  
\mathchardef\zh="7111  
\mathchardef\zvy="7112 
\mathchardef\zi="7113  
\mathchardef\zk="7114  
\mathchardef\zl="7115  
\mathchardef\zm="7116  
\mathchardef\zn="7117  
\mathchardef\zx="7118  
\mathchardef\zp="7119  
\mathchardef\zr="711A  
\mathchardef\zs="711B  
\mathchardef\zt="711C  
\mathchardef\zu="711D  
\mathchardef\zvf="711E 
\mathchardef\zq="711F  
\mathchardef\zc="7120  
\mathchardef\zw="7121  
\mathchardef\ze="7122  
\mathchardef\zy="7123  
\mathchardef\zf="7124  
\mathchardef\zvr="7125 
\mathchardef\zvs="7126 
\mathchardef\zf="7127  
\mathchardef\zG="7000  
\mathchardef\zD="7001  
\mathchardef\zY="7002  
\mathchardef\zL="7003  
\mathchardef\zX="7004  
\mathchardef\zP="7005  
\mathchardef\zS="7006  
\mathchardef\zU="7007  
\mathchardef\zF="7008  
\mathchardef\zW="700A  

\newcommand{\be}{\begin{equation}}
\newcommand{\ee}{\end{equation}}
\newcommand{\ra}{\rightarrow}
\newcommand{\lra}{\longrightarrow}
\newcommand{\bea}{\begin{eqnarray}}
\newcommand{\eea}{\end{eqnarray}}
\newcommand{\beas}{\begin{eqnarray*}}
\newcommand{\eeas}{\end{eqnarray*}}
\def\*{{\textstyle *}}
\newcommand{\R}{{\mathbb R}}
\newcommand{\T}{{\mathbb T}}
\newcommand{\C}{{\mathbb C}}
\newcommand{\unit}{{\mathbf 1}}
\newcommand{\SL}{SL(2,\C)}
\newcommand{\Sl}{sl(2,\C)}
\newcommand{\SU}{SU(2)}
\newcommand{\su}{su(2)}
\def\ssT{\sss T}
\newcommand{\G}{{\goth g}}
\newcommand{\D}{{\rm d}}
\newcommand{\Df}{{\rm d}^\zF}
\newcommand{\de}{\,{\stackrel{\rm def}{=}}\,}
\newcommand{\we}{\wedge}
\newcommand{\nn}{\nonumber}
\newcommand{\ot}{\otimes}
\newcommand{\s}{{\textstyle *}}
\newcommand{\ts}{T^\s}
\newcommand{\oX}{\stackrel{o}{X}}
\newcommand{\oD}{\stackrel{o}{D}}
\newcommand{\obD}{\stackrel{o}{\bD}}
\newcommand{\pa}{\partial}
\newcommand{\ti}{\times}
\newcommand{\A}{{\cal A}}
\newcommand{\Li}{{\cal L}}
\newcommand{\ka}{\mathbb{K}}
\newcommand{\find}{\mid}
\newcommand{\ad}{{\rm ad}}
\newcommand{\rS}{]^{SN}}
\newcommand{\rb}{\}_P}
\newcommand{\p}{{\sf P}}
\newcommand{\h}{{\sf H}}
\newcommand{\X}{{\cal X}}
\newcommand{\I}{\,{\rm i}\,}
\newcommand{\rB}{]_P}
\newcommand{\Ll}{{\pounds}}
\def\lna{\lbrack\! \lbrack}
\def\rna{\rbrack\! \rbrack}
\def\rnaf{\rbrack\! \rbrack_\zF}
\def\rnah{\rbrack\! \rbrack\,\hat{}}
\def\lbo{{\lbrack\!\!\lbrack}}
\def\rbo{{\rbrack\!\!\rbrack}}
\def\lan{\langle}
\def\ran{\rangle}
\def\zT{{\cal T}}
\def\tU{\tilde U}
\def\ati{{\stackrel{a}{\times}}}
\def\sti{{\stackrel{sv}{\times}}}
\def\aot{{\stackrel{a}{\ot}}}
\def\sati{{\stackrel{sa}{\times}}}
\def\saop{{\stackrel{sa}{\op}}}
\def\bwa{{\stackrel{a}{\bigwedge}}}
\def\svop{{\stackrel{sv}{\oplus}}}
\def\saot{{\stackrel{sa}{\otimes}}}
\def\cti{{\stackrel{cv}{\times}}}
\def\cop{{\stackrel{cv}{\oplus}}}
\def\dra{{\stackrel{\xd}{\ra}}}
\def\bdra{{\stackrel{\bd}{\ra}}}
\def\bAff{\mathbf{Aff}}
\def\Aff{\sss{Aff}}
\def\bHom{\mathbf{Hom}}
\def\Hom{\sss{Hom}}
\def\bt{{\boxtimes}}
\def\sot{{\stackrel{sa}{\ot}}}
\def\bp{{\boxplus}}
\def\op{\oplus}
\def\bwak{{\stackrel{a}{\bigwedge}\!{}^k}}
\def\aop{{\stackrel{a}{\oplus}}}
\def\ix{\operatorname{i}}
\def\V{{\cal V}}
\def\cD{{\cal D}}
\def\cC{{\cal C}}
\def\cE{{\cal E}}
\def\cL{{\cal L}}
\def\cN{{\cal N}}
\def\cR{{\cal R}}
\def\cJ{{\cal J}}
\def\cT{{\cal T}}
\def\cH{{\cal H}}
\def\bA{\mathbf{A}}
\def\bI{\mathbf{I}}
\def\wh{\widehat}
\def\wt{\widetilde}
\def\ol{\overline}
\def\ul{\underline}
\def\Sec{\sss{Sec}}
\def\Lin{\sss{Lin}}
\def\ader{\sss{ADer}}
\def\ado{\sss{ADO^1}}
\def\adoo{\sss{ADO^0}}
\def\AS{\sss{AS}}
\def\bAS{\sss{AS}}
\def\bLS{\sss{LS}}
\def\bAP{\sss{AV}}
\def\bLP{\sss{LP}}
\def\AP{\sss{AP}}
\def\LP{\sss{LP}}
\def\LS{\sss{LS}}
\def\Z{\mathbf{Z}}
\def\oZ{\overline{\bZ}}
\def\oA{\overline{\bA}}
\def\cim{{C^\infty(M)}}
\def\de{{\cal D}^1}
\def\la{\langle}
\def\ran{\rangle}
\def\by{{\bi y}}
\def\bs{{\bi s}}
\def\bc{{\bi c}}
\def\bd{{\bi d}}
\def\bh{{\bi h}}
\def\bD{{\bi D}}
\def\bY{{\bi Y}}
\def\bX{{\bi X}}
\def\bL{{\bi L}}
\def\bV{{\bi V}}
\def\bW{{\bi W}}
\def\bS{{\bi S}}
\def\bT{{\bi T}}
\def\bC{{\bi C}}
\def\bE{{\bi E}}
\def\bF{{\bi F}}
\def\bP{{\bi P}}
\def\bp{{\bi p}}
\def\bz{{\bi z}}
\def\bZ{{\bi Z}}
\def\bq{{\bi q}}
\def\bQ{{\bi Q}}
\def\bx{{\bi x}}

\def\sA{{\sss A}}
\def\sC{{\sss C}}
\def\sD{{\sss D}}
\def\sG{{\sss G}}
\def\sH{{\sss H}}
\def\sI{{\sss I}}
\def\sJ{{\sss J}}
\def\sK{{\sss K}}
\def\sL{{\sss L}}
\def\sO{{\sss O}}
\def\sP{{\sss P}}
\def\sPh{{\sss P\sss h}}
\def\sT{{\sss T}}
\def\sV{{\sss V}}
\def\sR{{\sss R}}
\def\sS{{\sss S}}
\def\sE{{\sss E}}
\def\sF{{\sss F}}
\def\st{{\sss t}}
\def\sg{{\sss g}}
\def\sx{{\sss x}}
\def\sv{{\sss v}}
\def\sw{{\sss w}}
\def\sQ{{\sss Q}}
\def\sj{{\sss j}}
\def\sq{{\sss q}}
\def\xa{\tx{a}}
\def\xc{\tx{c}}
\def\xd{\tx{d}}
\def\xD{\tx{D}}
\def\xV{\tx{V}}
\def\xF{\tx{F}}
\def\dt{\xd_{\sss T}}
\def\vt{\textsf{v}_{\sss T}}
\def\vta{\operatorname{v}_\zt}
\def\vtb{\operatorname{v}_\zp}
\def\cM{\cal M}
\def\cN{\cal N}
\def\cD{\cal D}
\def\ug{\ul{\zg}}
\def\sTn{\stackrel{\scriptscriptstyle n}{\textstyle\sT}\!}
\def\sTd{\stackrel{\scriptscriptstyle 2}{\textstyle\sT}\!}
\def\stn{\stackrel{\scriptscriptstyle n}{\textstyle\st}\!}
\def\std{\stackrel{\scriptscriptstyle 2}{\textstyle\st}\!}
\newdir{ (}{{}*!/-5pt/@^{(}}


\setcounter{page}{1} \thispagestyle{empty}


\bigskip

\bigskip

\title{The Tulczyjew triple for classical fields
}

        \author{
        Katarzyna  Grabowska\thanks{The research financed by the Polish Ministry of Science and Higher Education under the
 grant N N201 365636.}\\ \\
          {\it Faculty of Physics}\\
                {\it University of Warsaw} \\
                {Ho\.za 69, 00-681 Warszawa, Poland}}

\date{}
\maketitle
\begin{abstract}
{The geometrical structure known as the Tulczyjew triple has proved to be very useful in describing mechanical
systems, even those with singular Lagrangians or subject to constraints. Starting from basic concepts of
variational calculus, we construct the Tulczyjew triple for first-order Field Theory. The important feature of
our approach is that we do not postulate {\it ad hoc} the ingredients of the theory, but obtain them as
unavoidable consequences of the variational calculus. This picture of Field Theory is covariant and complete, containing not
only the Lagrangian formalism and Euler-Lagrange equations but also the phase space, the phase dynamics and
the Hamiltonian formalism. Since the configuration space turns out to be an affine bundle, we have to use
affine geometry, in particular the notion of the affine duality. In our formulation,
the two maps $\alpha$ and $\beta$ which constitute the Tulczyjew triple are morphisms of double structures of
affine-vector bundles. We discuss also the Legendre transformation, i.e. the transition between the Lagrangian
and the Hamiltonian formulation of the first-order field theory.}


\end{abstract}

\section{Introduction}
Variational calculus is a natural language for describing statics of mechanical systems. All mathematical
objects that are used in statics have { direct} physical interpretations. Moreover, similar mathematical tools
are widely used also in other theories, like dynamics of particles or field theories. In classical mechanics
variational calculus was used first for deriving equations of motion of mechanical system, i.e. the
Euler-Lagrange equations.

In numerous works by W. M. Tulczyjew, for example in the book \cite{TU3} and papers \cite{Tu5, Tu6, Tu7, Tu9,
Tu10}, one may find another philosophy of using variational calculus in mechanics and field theories. This
philosophy, especially { the one leading to the construction} called the Tulczyjew triple, has been recently
recognized by many theoretical physicists and mathematicians. The main advantage of the approach developed by
Tulczyjew and his collaborators is its generality. For example, using the Tulczyjew triple for autonomous
mechanics we can derive the phase equations for systems with singular Lagrangians and understand { properly}
the Hamiltonian description of such systems. One can even discuss systems with more general generating objects
than just a Lagrangian function, { e.g. systems} described by family of Lagrangians or a Lagrangian function
defined on a submanifold. For the details we refer to \cite{TU}.

Another advantage of the { Tulczyjew's} approach is { its} flexibility. { Being based on well-defined general
principles, it can be easily adapted to different settings. No wonder that there} were many attempts to
generalize the Tulczyjew triple to more general contexts of mechanics on algebroids or different field
theories (see e.g \cite{EM,LMS,RRR}). In our earlier paper \cite{G} we started with constructing a toy model
of the triple for field theory in the simplest topological situation.

The purpose of this work is to construct the Tulczyjew triple for first-order field theory in a very general
setting, i.e. in the case where fields are sections of some differential fibration with no additional
structure assumed.  The origins of the geometric structures we study lie in the rigorous formulation of the variational
principle including boundary terms. We pay much attention to recognize physically important objects, like the
phase space, phase dynamics, the Legendre map, Hamiltonians, etc. These issues are usually not well elaborated
in the literature, as the classical field theory models use to concentrate on the Euler-Lagrange equations. Of
course, we recover also the commonly accepted Euler-Lagrange equations, this time without requiring any
regularity of the Lagrangian.

Classical field theory is usually associated to the concept of multisymplectic structure.  The literature on
the subject is very rich, so we mention only a few main papers. The multisymplectic approach appeared
first in \cite{Tu5} and \cite{Kij2}. Then, it was developed by Gotay, Isennberg, Marsden and others in \cite{GIM1,
GIM2, G1, G2}. The original idea of the multisymplectic structure has been thoroughly investigated and
developed by many authors, see { e.g.} papers by Cari\~nena, Crampin, Ibort, Cantrijn, De Leon \cite{CIL,CCI1,
CCI2} and Echeverria-Enriquez, Mu\~noz-Lecanda \cite{EM} for general analysis of the multisymplectic structure
and its application to the classical field theories, and by Forger, Paufler, R\"{o}mer \cite{FP1, FP2}, or
Vankershaver, Cantrijn, De Leon \cite{VCL} for the discussion of more detailed problems associated to the
structure. An interesting discussion of the problem can be found also in the paper \cite{HK}. The Tulczyjew
triple in the context of multisymplectic field theories appeared recently in \cite{LMS}. A similar diagram, however
with differences on the Hamiltonian side, one can find also in \cite{GM}. Another approach to field theory,
based on differential forms on fibre bundles, is present in works by Krupkov\'{a} and collaborators,
e.g. \cite{Kr1, Kr2,Kr3}).

Our approach to the Lagrangian and Hamiltonian formalism developed in the paper is different. We do not use
directly the multisymplectic formalism, building instead the triple out of natural morphisms of double
structures of affine-vector bundles.  Also, we do not use the framework based on Klein's ideology and do not
concentrate on the Euler-Lagrange equations nor regular Lagrangians, since the phase dynamics is for us the
principal object. Using the affine geometry as a tool and following guide-lines of variational calculus we
arrive to spaces and maps on the Hamiltonian side of the triple. The variational problem we start with determines
uniquely the phase space together with its canonical structure which is different from the one in \cite{GM}. Moreover, the
canonical structure of the phase space is not a multisymplectic form, but a family of symplectic forms on
fibres over the base manifold with values in the space of forms on the base. The two structures are related,
but not identical. As far as we know, a similar research is being done independently, for instance, by L.
Vitagliano \cite{V}, and E. Guzm\'an. Just before submitting this paper we spotted a preprint by Campos,
Guzm\'an, and Marrero \cite{CGM} dealing with similar questions.

The starting point of our studies is a locally trivial fibration $\zz:E\ra M$ over a manifold $M$ of dimension $m$, whose
sections represent fields, and the corresponding bundle $J^1E$ of first jets of sections playing the role of
kinematic configurations. A Lagrangian is a map $L:\sJ^1E\to \Omega^m$, where $\Omega^k:=\bigwedge^k\sT^*M$ is
the bundle of $k$-forms on $M$. The phase space of the theory turns out to be the bundle
$\mathcal{P}=\sV^*E\ot_E\zz^*(\zW^{m-1})$ denoted simply $\mathcal{P}=\sV^*E\ot_E\zW^{m-1}$. Here, $\sV^*E$ is
the dual of the vertical subbundle $\sV E$ in $\sT E$ and $\zz^*(\zW^{m-1})$ is the pull-back bundle of
$\zW^{m-1}$ along the projection $\zz:E\ra M$.

The Lagrangian side of the Tulczyjew triple is constituted by a map
$$\za:\sJ^1\mathcal{P}\to \sV^*\sJ^1E\ot_{\sJ^1E}\zW^{m}$$
being a morphism of double structures of affine-vector bundles associated with fibrations over $\mathcal{P}$
and $\sJ^1 E$. The vertical derivative $\xd L$ of the Lagrangian is a section of the bundle
$\sV^*\sJ^1E\ot_{\sJ^1E}\zW^{m}\to \sJ^1E$ and the (implicit) phase dynamics is defined as a subset
$\mathcal{D}$ of $\sJ^1\mathcal{P}$ being the inverse image by $\za$ of the image of $\xd L$, i.e.,
$\mathcal{D}=\za^{-1}(\xd L(\sJ^1E))$.

Similarly, the Hamiltonian side of the triple is built on another  morphism of affine-vector bundles, fibred
also over $\mathcal{P}$ and $\sJ^1 E$,
$$\zb:\sJ^1\mathcal{P}\to \sP\sJ^\dag E\,,$$
where $\sJ^\dag E$ is the `affine dual' of $\sJ^1E$, i.e. the bundle of affine maps from $\sJ^1_eE$ to
$\zW^m_{\zz(e)}$, and $\sP\sJ^\dag E$ is the affine phase bundle of the affine line bundle
$\zvy:\sJ^1E\to \mathcal{P}$, an affine analog of the cotangent bundle. A Hamiltonian is a section of the
bundle $\zvy$, i.e. $H:\mathcal{P}\to \sJ^1E$, so its affine differential $\xd H$ can be viewed as a map $\xd
H:\mathcal{P}\to\sP\sJ^\dag E$ and defines a phase dynamics $\mathcal{D}=\zb^{-1}(\xd H(\mathcal{P}))$. It is
obvious that, being presented in a coordinate-free form, the whole theory is covariant.

We would like to stress that all the geometrical objects we construct are not just postulated {\it ad hoc},
but discovered by starting from natural general principles and rigorous investigations of the geometry which
arises in this way. Another nice feature of our approach is that it admits a straightforward generalization
for field theories of higher orders. This issues, however, we postpone to a separate paper.

This work is organized as follows. In section \ref{secintro} we present the conceptual  background of
variational theories which justifies mathematical constructions and physical interpretations in particular
examples. After introducing some notation in section \ref{sec2}, we pass to the Lagrangian side of the triple in section
\ref{seclagr}. Section \ref{secham} is devoted to the Hamiltonian side of the triple.

\section{Variational calculus in Statics and Mechanics}\label{secintro}

Let us start with the simplest case of statics of a mechanical system. We shall assume that the set of all
possible configurations of the system is a differential manifold $Q$. The tangent and cotangent bundles
$$\tau_Q: \sT Q\longrightarrow Q,\qquad\pi_Q: \sT^\ast Q\longrightarrow Q$$
will also be used. In statics we are usually interested in equilibrium configurations of an isolated system,
as well as a system with an interaction with other static systems. The system alone or in interaction is
examined by preforming processes and calculating the cost of each process. We assume that all the processes
are {\it quasistatic}, i.e. they are slow enough to produce no dynamical effects. Every process can be
represented by a one-dimensional smooth oriented submanifold with boundary. It may happen that for some
reasons, not all the processes are admissible, i.e. the system is constrained. All the information about the
system is therefore encoded in three objects: the configuration manifold $Q$, the set of all admissible
processes, and the cost function that assigns a real number to every process. The cost function should fulfill
some additional conditions, e.g. it should be additive in the sense that if we break a process into two
subprocesses, then the cost of the whole process should be equal to the sum of the costs of the two
subprocesses. Usually we assume that the cost function is local, i.e. for each process it is an integral of a
certain positively homogeneous function $W$ on $\sT Q$. There are distinguished systems, called {\it regular},
for which all the processes are admissible and the function $W$ is the differential of a certain function
$U:Q\rightarrow \R$. In this case $U$ is called the {\it internal energy} function.

An {\it equilibrium point} for the system is such a point $q\in Q$ that all the processes starting from $q$
have positive cost, at least initially, i.e. for some sufficiently small subprocess with the same initial
point. Usually we formulate only a first-order necessary criterion for the equilibrium point. It says that a
point $q$ is an equilibrium point of the system if
$$W(\delta q)\geq 0$$
for all vectors $\delta q\in\sT_q Q$ tangent to admissible processes. Vectors tangent to admissible processes
are called {\it admissible virtual displacements}. The set of such vectors will be denoted with $\Delta$. It
may happen that $\Delta\cap\sT Q$ does not project on the whole $Q$. We have then the set of {\it admissible
configurations} $C=\tau_Q(\Delta\cap\sT Q)$. For regular systems the equilibrium condition assumes the form
$$\xd U(q)=0.$$

We examine the interaction between two systems by creating composed systems. We can compose systems that have the
same configuration space $Q$. The composite system is described by the intersection of the sets of admissible
processes and the sum of the cost functions. We describe our system by making a list of all systems that, composed
with our system, have certain admissible configuration $q$ as an equilibrium point. We observe that at each
admissible $q$ all the external systems interacting with our systems can be classified according to their
influence on our system. Moreover, in every class we can find a regular system, therefore the whole
class can be represented by the differential of the internal energy of that regular system. We call {\it a
force} the class of external systems interacting with our system. The force is represented by a covector, i.e.
an element of $\sT^\ast Q$. Instead of making a list of all external systems in equilibrium with our system at
the point $q$, we can give a subset of $\sT_q^\ast Q$ representing those systems. The subset of $\sT^\ast Q$
of all forces in equilibrium with our system at all admissible points we call the {\it constitutive set}. For
a large class of static systems the constitutive set contains the complete information about the system. The
passage from the triple $(Q, \Delta\subset\sT Q, W)$ describing our system to the constitutive set $D\subset
\sT^\ast Q$ is called the {\it Fenchel-Legendre transformation}.

Let us now use the above concepts to describe the autonomous dynamics of a non-relativistic particle. For
simplicity, we will consider only the unconstrained case. Let us assume that the set of positions of the
particle (possible configurations) is a smooth manifold $Q$. There are at least two approaches to the problem.
The first deals with the finite time interval $[t_0, t_1]$, while the second with the infinitesimal time
interval represented by the Dirac $\delta$-distribution at $t$. The finite case provides a useful
representation of objects coming from statics, while infinitesimal approach leads to differential equations
for phase trajectories that are commonly used in physics. We skip the details of the construction and provide
here only a summary of the results of both approaches.

For the finite time interval the configuration space is the space of all motions, i.e. pieces of smooth curves
parameterized by the time, $\gamma: [t_0, t_1]\ni t\rightarrow Q$. This space is not a standard manifold any
more, therefore we have to precise the notions of a smooth function, a tangent vector and a covector. We will work
with functions (usually called functionals) of the form
\be S(\gamma)=\int_{t_0}^{t_1}L(\st\gamma(t))\xd t.\label{intro1}\ee
In the above formula, $L$ is a smooth function on $\sT Q$, called the {\it Lagrangian}, and $\st \gamma$
denotes the tangent prolongation of the curve $\gamma$. A curve in the configuration space always comes from a
homotopy, i.e. from a smooth map
$$\chi:\R^2\ni(s,t)\longmapsto \chi(s,t)\in Q.$$
Restricting the domain of $t$ to $[t_0, t_1]$ for every $s$, we obtain a curve in the space of motions that is
smooth by definition. The choice is justified by the fact that the composition of the curve with  any function
of the form (\ref{intro1}) is a real function smooth in the usual sense.

A vector tangent to a manifold is usually defined as an equivalence class of curves. In our situation we can
adopt the same definition. Working with equivalence classes is difficult, therefore we observe that each
equivalence class at a configuration $\gamma$ can be conveniently represented as a curve
\be\delta\gamma: [t_0, t_1]\longrightarrow\sT Q\ee
such that $\tau_Q\circ\delta\gamma=\gamma$. In differential geometry we define covectors as equivalence
classes of functions. Equivalence classes are again too abstract objects, therefore we need a convenient
representation for covectors. The idea of such a representation is given by performing variation of the
functional $S$ and separating boundary terms like in the procedure of deriving the Euler-Lagrange equations:
\be\langle\delta S,\delta \gamma\rangle=
\int_{t_0}^{t_1}\langle\mathcal{E}L(\st^2\gamma(t)),\delta\gamma(t)\rangle\xd t+
\langle\mathcal{P}L(\st\gamma(t_1), \delta\gamma(t_1)\rangle -\langle\mathcal{P}L(\st\gamma(t_0),
\delta\gamma(t_0)\rangle,\ee
where $\mathcal{E}L$ denotes the Euler-Lagrange variation of $L$ that depends on
the second prolongation $\st^2\gamma$ of the motion $\gamma$, and $\mathcal{P}L$ is a vertical differential of
$L$ with respect to the projection $\tau_Q$. We see that the convenient representation of a covector is a
triple $(f,p_0,p_1)$, where
\be f:[t_0,t_1]\rightarrow \sT^\ast Q,\qquad \pi_Q\circ f=\gamma,\qquad p_0\in\sT_{\gamma(t_0)}Q,\qquad
p_1\in\sT_{\gamma(t_1)}Q.\ee The evaluation between $\delta\gamma$ and $(f,p_0, p_1)$ is given by
\be\langle\!\langle(f, p_0, p_1),\delta\gamma\rangle\!\rangle=\int_{t_0}^{t_1}\langle f(t),\delta\gamma(t)\rangle\xd t+
\langle p_1, \delta\gamma(t_1)\rangle -\langle p_0, \delta\gamma(t_0)\rangle.\ee
The symbol
$\langle\!\langle\cdot, \cdot\rangle\!\rangle$ will denote the pairing between convenient representations.

The elements $f$, $p_0$, and $p_1$ have physical interpretation. The curve $f$ is an external force acting on
the particle during its motion, $p_0$ and $p_1$ are the initial and the final momenta. The constitutive set
consists of all triples $(f,p_0,p_1)$ such that the particle moves along the curve $\gamma=\pi_Q\circ f$,
starting at $\gamma(t_0)$ with the initial momentum $p_0$ and arriving to $\gamma(t_1)$ with the final momentum $p_1$,
while it is a subject to the force $f$ along the motion. Having the constitutive set, we can discuss isolated
system, i.e. systems with the external force $f=0$, as well as the system interacting with external forces of
different kinds. We may interpret the momenta as a result of an interaction between the system and its past
and its future. Therefore it makes no sense to keep the momenta equal to zero. The space of momenta is usually
called the {\it phase space} of the system. We see that in our example the phase space is $\sT^\ast Q$.

The constitutive set, as defined above, is a complicated object. We would like to describe it in a more
convenient way, e.g. using differential equations for curves in forces and momenta such that their solutions
restricted to any interval $[t_0, t_1]$ lie in the constitutive set. We can obtain such equations using the infinitesimal
approach to the dynamics.

Passing to the infinitesimal formulation, we replace the finite domain of the integration, $[t_0, t_1]$, with the
Dirac's $\delta$-distribution at the point $t$. We see that the configurations are now elements of $\sT Q$. Since
the configuration space is again a manifold, we have natural notions of smooth functions, curves, tangent
vectors, and covectors. The "internal energy" function is now just the Lagrangian, and its differential is an
element of $\sT^\ast\sT Q$. Virtual displacements are vectors tangent to curves in $\sT Q$, i.e. elements of
$\sT\sT Q$. We observe, however, that studying convenient representations of vectors and covectors for the finite
formulation gives interesting results also in the infinitesimal limit. A virtual displacement of a configuration
$\st\gamma(t)$ is a vector $\delta\st\gamma(t)$ in $\sT\sT Q$ such that
\be \tau_{\sT Q}(\delta\st\gamma(t))=\st\gamma(t), \qquad\sT\tau_Q(\delta\st\gamma(t))=\delta\gamma(t).\ee
In turn, the infinitesimal limit of the virtual displacement $\delta\gamma$ that we had for the finite time interval
is an element $\st\delta\gamma(t)\in \sT\sT Q$ such that
\be \tau_{\sT Q}(\st\delta\gamma(t))=\delta\gamma(t), \qquad\sT\tau_Q(\st\delta\gamma(t))=\st\gamma(t).\ee
The virtual displacement and its convenient representation are two elements of $\sT\sT Q$ related by the
canonical flip
$$\kappa_M: \sT\sT Q\longrightarrow \sT\sT Q.$$

The constitutive set for our regular system in infinitesimal setting is the graph of $\xd L$. Again, the
convenient representation of elements of the cotangent bundle in the finite time interval formulation provides
us with another useful interpretation of the constitutive set. For the finite time interval the evaluation of
$(f, p_0, p_1)$ on $\delta\gamma$ reads
\be \int_{t_0}^{t_1}\langle f(t), \delta\gamma\rangle \xd t+
\langle p_1, \delta\gamma(t_1)\rangle- \langle p_0, \delta\gamma(t_0)\rangle.\ee
In the infinitesimal case we get
\be\langle f(t), \delta\gamma(t)\rangle+\frac{\xd}{\xd t}\langle p(t), \delta\gamma(t)\rangle.\ee
In the absence of external forces, another description of a constitutive set can be derived out of the equation
\be\frac{\xd}{\xd t}\langle p(t), \delta\gamma(t)\rangle=\langle\xd L(\st\gamma(t)), \delta\st\gamma(t)\rangle\label{intro2}\ee
The left-hand side is the so called {\it tangent evaluation} between a vector tangent to $\sT^\ast Q$ and a vector
tangent to $\sT Q$ such that they have common tangent projection on $\sT Q$. More precisely, if $p:
\R\rightarrow \sT^\ast Q$ and $\delta\gamma:\R\rightarrow\sT Q$ are two curves covering the same curve
$\gamma:\R\rightarrow Q$, then
\be\langle\!\langle \st p(t), \st\delta\gamma(t)\rangle\!\rangle=
\frac{\xd}{\xd t}\langle p(t), \delta\gamma(t)\rangle.\ee
Since the vectors $\st\delta\gamma(t)$  and
$\delta\st\gamma(t)$ are related by the canonical flip $\kappa_Q$,
$$\kappa_Q(\delta\st\gamma(t))=\st\delta\gamma(t),$$
the differential of the Lagrangian $\xd L(\st\gamma(t))$ and the tangent vector $\st p(t)$ are related by the
Tulczyjew $\alpha_Q$ (which is dual to $\kappa_Q$),
$$\alpha_Q:\sT\sT^\ast Q\longrightarrow \sT^\ast\sT Q.$$
In this way we have obtained another description of the constitutive set, called {\it the phase dynamics} and
given by the formula
$$\sT\sT^\ast M\supset D=\alpha_Q^{-1}(\xd L(\sT Q)).$$

If the system is autonomous, then the constitutive set for any time $t$ is the same. The condition for a curve
$p: \R\supset I\rightarrow \sT^\ast Q$ to be a phase trajectory of the system is that
$$\forall t\in I\quad\st p(t)\in D.$$
The dynamics $D$ can be understood as a differential-algebraic equation for a pair of curves, $f:
\R\rightarrow \sT^\ast Q$ and $p: \R\rightarrow \sT^\ast Q$, covering the same curve in $Q$.
A curve $\gamma: I\rightarrow Q$ satisfies, in turn, the corresponding {\it Euler-Lagrange equation}, if the curve
$I\ni t\mapsto\alpha_Q^{-1}(\xd L(\zg(t)))\in\sT\sT^\ast Q$ is the tangent prolongation of its projection to
$\sT^\ast Q$.

External forces can be included in the picture as follows. Equation (\ref{intro2}) completed with the force
reads as
\be\langle f(t), \delta\gamma(t)\rangle+
\frac{\xd}{\xd t}\langle p(t), \delta\gamma(t)\rangle=\langle\xd L(\st\gamma(t)),
\delta\st\gamma(t)\rangle\label{intro3}.\ee
The set of all elements of $\sT\sT^\ast Q$ with fixed projections
on $\sT^\ast Q$ and $\sT Q$ is an affine space modeled on the appropriate fibre of $\sT^\ast Q$. The force $f$
can be therefore added to a vector tangent to the phase space. The map $\alpha_Q$ can now be extended to
the map
$$\widetilde\alpha_Q:\sT^\ast Q\times_Q\sT\sT^\ast Q\longrightarrow \sT^\ast\sT Q, \qquad
\widetilde\alpha_Q(f,u)=\alpha_Q(u+f).$$ The dynamics with external forces (see \cite{MTU}) is a subset
$\widetilde D$ of $\sT^\ast Q\times_Q\sT\sT^\ast Q$ given by
$$\widetilde D=\widetilde\alpha_Q{}^{-1}(\xd L(\sT Q)).$$

All the structures needed for generating the dynamics from a Lagrangian (without external forces) can be
summarized in the following commutative diagram of vector bundle morphisms
\be\label{intro5}\xymatrix@C-5pt{
& \sT\sT^\ast Q \ar[rrr]^{\alpha_Q} \ar[rdd]^{\sT\pi_Q} \ar[ld]_{\zt_{\sT^\ast Q}} & & & \sT^\ast\sT Q
\ar[rdd]^{\pi_{\sT Q}} \ar[ld]_{\zeta}\\
\sT^\ast Q\ar[rrr]^/-25pt/{id}\ar[rdd]^/-20pt/{\pi_Q} & & & \sT^\ast Q\ar[rdd]^/-20pt/{\pi_Q} &  \\
 & & \sT Q\ar[rrr]^/-25pt/{id}\ar[ld]_{\zt_Q} & & & \sT Q\ar[ld]_{\zt_Q} \\
&  Q \ar[rrr]^{id}& & &  Q }\ee The map $\alpha_Q$ is an isomorphism of double vector bundles. Recall that
double vector bundles are, roughly speaking, manifolds equipped with two compatible vector bundle structures.
The compatibility condition can be expressed e.g. as the commutation of the two Euler vector fields associated with
these vector bundle structures. A precise definition of a double vector bundle together with its basic properties
can be found in \cite{Pr1,KU,GR}.

The map $\alpha_Q$ is also a symplectomorphism between the symplectic manifolds $(\sT\sT^\ast Q,
\xd_\sT\omega_{Q})$ and $(\sT^\ast\sT Q, \omega_{\sT Q})$, where $\xd_\sT\omega_{Q}$ is the complete lift of
the canonical symplectic form $\omega_Q$ on $\sT^\ast Q$ and $\omega_{\sT Q}$ is the canonical symplectic form
on $\sT^\ast\sT Q$.

It may happen that the phase dynamics is an implicit differential equation, i.e. it is not the image of a
vector field. In some cases, however, the phase dynamics is the image of a Hamiltonian vector field for some
function $H:\sT^\ast Q\rightarrow \R$. So that we can write
\be\mathcal{D}=\beta_Q^{-1}(\xd H(\sT^\ast Q)),\ee
where $\beta_Q$ is the canonical isomorphism between $\sT\sT^\ast Q$ and $\sT^\ast\sT^\ast Q$ given by the
canonical symplectic form $\omega_Q$ on $\sT^\ast Q$,
$$\beta_Q:\sT\sT^\ast Q\longrightarrow\sT^\ast\sT^\ast Q, \qquad\langle\beta_Q(v),w\rangle=\omega_Q(v,w).$$
Let us recall that the canonical symplectic form $\omega_Q$ is defined as the differential
\begin{equation}\label{intro6}\omega_Q=\xd \vartheta_Q\end{equation}
of the Liouville form $\vartheta_Q$ given by
\begin{equation}\label{intro7}\vartheta_Q(v)=\langle \zt_{\sT^\ast Q}(v), \sT\pi_Q(v) \rangle.\end{equation}

The structures needed for Hamiltonian mechanics can be presented in the following commutative diagram:
\be\label{intro8}\xymatrix@C-5pt{
& \sT^\ast\sT^\ast Q  \ar[rdd]^{\zx} \ar[ld]_{\pi_{\sT^\ast Q}} & & & \sT\sT^\ast Q\ar[lll]_{\beta_Q}
\ar[rdd]^{\sT\pi_{Q}} \ar[ld]_{\zt_{\sT^\ast Q}}\\
\sT^\ast Q\ar[rdd]^/-20pt/{\pi_Q} & & & \sT^\ast Q\ar[rdd]^/-20pt/{\pi_Q}\ar[lll]_/-25pt/{id} &  \\
 & & \sT Q\ar[ld]_{\zt_Q} & & & \sT Q\ar[ld]_{\zt_Q}\ar[lll]_/-25pt/{id} \\
&  Q & & &  Q\ar[lll]_{id} }\ee The map $\beta_Q$ is an isomorphism of double vector bundles.

The formulation of the autonomous mechanics described above has at least two important features when compared
with the ones in textbooks: it is very simple and can be easily generalized to more complicated cases
including constraints, nonautonomous mechanics, and mechanics on algebroids \cite{GGU0,GGU1, GGU2}. And last but
not least, we need no regularity conditions for the Lagrangian. The Lagrangian  and Hamiltonian can be
functions, but one of them (or both) can be replaced by families of functions generating  Lagrangian
submanifolds in $\sT^\ast\sT Q$ and $\sT^\ast\sT^\ast Q$, respectively. It happens e.g. in the case of a
relativistic particle in the Minkowski space \cite{TU}. Moreover, the generating object for dynamics on the
Lagrangian side can be replaced by a $1$-form different from $\xd L$. It happens
e.g. for systems with friction. The crucial role is played by two mappings: $\alpha_Q$ and $\beta_Q$.

The two diagrams (\ref{intro7}) and (\ref{intro8}) glued together are called the {\it Tulczyjew triple} for
mechanics. We would like to emphasize the fact that the Lagrangian side of the Tulczyjew triple is not
postulated, but derived form the variational calculus. The Hamiltonian side, present only in infinitesimal
formulation, comes from the fact that $\sT\sT^\ast Q$ is equipped with two Liouville  structures, i.e.  is
isomorphic to two different cotangent bundles $\sT^\ast\sT^\ast Q$ and $\sT^\ast\sT Q$ (see \cite{TU2}).

The paper is devoted to deriving the Tulczyjew triple for field theory, i.e. to the case where configurations
are sections of a certain fibration. Like in mechanics, the Lagrangian side of the triple appears as a result
of the existence of  so called convenient representations of equivalence classes representing tangent vectors
and covectors. The Hamiltonian side is again related to the canonical isomorphism between certain cotangent
bundles.

\section{Notation}\label{sec2}
The notation used in papers concerning the geometry of classical field theory is usually very complicated,
because iterated tangent functors have to be used. In this section we will present a system of notation that
will be used in the following sections. We will try to introduce some rules which make the notation systematic for
the cost of the length of some symbols.
\smallskip

For $M$ being a smooth manifold of dimension $m$, we denote by
$$\zt_M: \sT M\longrightarrow M\qquad\text{and}\qquad  \zp_M:\sT^\ast M\longrightarrow M$$
the tangent and cotangent bundles. Let $U\subset M$ be a domain of local coordinate system $(q^i)_{i=1}^m$ on
$M$. We have the adapted systems of coordinates $(q^i,\dot q^j)$ and $(q^i,p_j)$ on $\zt_M^{-1}(U)$ and
$\zp_M^{-1}(U)$, respectively. The manifold $M$ is the manifold on which the field is defined, in applications
it can be e.g. the space-time. We will assume in the following that the manifold $M$ is oriented. This
assumption allows us to use even forms instead of densities and odd forms (which are not commonly used).
However, we will relax this assumption in the example to have clear physical interpretations of geometrical
objects. Instead of the usual notation $\bigwedge^k\sT^\ast M$ for the space of $k$-covectors on $M$, we will
use the shorter symbol $\Omega^k$. The space $\Omega^k$ is a vector bundle over $M$. The vector space over the
point $x\in M$ will be denoted by $\Omega^k_x$. For a local coordinate system $(q^i)$ there is associated the
volume form $\eta=\xd q^1\wedge\xd q^2\wedge\cdots\wedge \xd q^m$. We will denote by $\eta_i$ the contraction
$\imath(\frac{\partial}{\partial q^i})\eta$.
\smallskip

Let
$$\zz: E\rightarrow M$$
be a smooth locally trivial fibration with the total space of dimension $m+n$. The total space is the space of
values of the field, i.e. a field is a local section of the bundle $\zz$. In applications the bundle $\zz$
can have additional structures. On an open subset $V\subset E$ such that $\zz(V)=U$ we can introduce a local
coordinate system $(q^i, y^a)$ adapted to the structure of the bundle.
\smallskip

The space of vectors tangent to $E$ and vertical with respect to the projection onto $M$ will be denoted by $\sV
E$. For the restriction of $\zt_E$ to the space of vertical vectors we will use the symbol $\zn_E$. The bundle
$$\zn_E: \sV E\longrightarrow E$$
is therefore a vector bundle. The adapted coordinate system on $\zn_E^{-1}(V)$, coming from $(q^i, y^a)$ on
$V$, will be denoted by $(q^i, y^a, \delta y^b)$. We will need also the dual vector bundle
$$\zr_E:\sV^\ast E\rightarrow E$$
with local coordinates $(q^i, y^a, p_b)$ defined on $\zr_E^{-1}(V)$.
\smallskip

The space of first jets of sections of the bundle $\zz$ will be denoted $\sJ^1 E$. By definition, the first
jet $\sj^1\zs(e)$ of the section $\zs$ at the point $e$ is an equivalence class of sections having the same
value $e$ at the point $x=\zz(e)$ and such that the spaces tangent to the graphs of the sections at the point
$e$ coincide. Therefore, there is a natural projection $\sj^1\zz$ from the space $\sJ^1\zz$ onto the manifold
$E$,
$$\sj^1\zz:\;\;\sJ^1 E\ni\sj^1\sigma(e)\longmapsto e\in E.$$
Moreover, every jet $\sj^1\sigma(e)$ can be identified with the linear map from $\sT_{x}M$ to $\sT_e E$ being
the tangent map $\sT\zs$ restricted to the space $\sT_x M$. It is easy to see that linear maps coming from
jets at the point $e$ form an affine subspace of the space $Lin(\sT_{x}M, \sT_e E)$ of all linear maps from
$\sT_x M$ to $\sT_e E$.
 A map belongs to this subspace if, composed with
$\sT\zz$, gives the identity. Using the tensorial representation, we can therefore write the inclusion
$$\sJ^1_e E\subset\sT^\ast_xM\otimes\sT_eE.$$
The affine subspace $\sJ^1_e E$ is modeled on a vector subspace of maps having values in vectors tangent to
$E$ at the point $e$ and vertical with respect to the fibration $\zz$. In tensor representation the model
vector space is $\sT^\ast_xM\otimes\sV_eE$. Summarizing, the bundle $\sJ^1 E \rightarrow E$ is an affine
subbundle in the vector bundle
$$\zz^\ast(\sT^\ast M)\otimes_E\sT E\longrightarrow E$$
modeled on the vector bundle
$$\zz^\ast(\sT^\ast M)\otimes_E\sV E\longrightarrow E.$$
The symbol $\zz^\ast(\sT^\ast M)$ denotes the pull-back of the cotangent bundle $\sT^\ast M$ along the
projection $\zz$. In the following, we will omit the symbol of the pull-back, writing simply $\sT^\ast
M\otimes_E\sT E$ and $\sT^\ast M\otimes_E\sV E$.

\smallskip

Using the adapted coordinates $(q^i, y^a)$ in $V\subset E$, we can construct the induced coordinate system
$(q^i, y^a, y^b{}_j)$ on $(\sJ^1\zz)^{-1}(V)$ such that, for any section $\sigma$ given by $n$ functions
$\sigma^a(q^i)$, we have
$$y^b{}_j(\sj^1\sigma(e))=\frac{\partial\sigma^b}{\partial q^j}(q^i(x)).$$
In the tensorial representation, the jet $\sj^1\sigma(e)$ can be written as
$$\xd q^i\otimes\frac{\partial}{\partial q^i}+
\frac{\partial\sigma^a}{\partial q^j}(q^i(x))\xd q^j\otimes\frac{\partial}{\partial y^a},$$ where we have used
local bases of sections of $\sT^\ast M$ and $\sT E$ coming from the chosen coordinates.
\smallskip

In the following, we will have to use iterated tangent functors $\sJ^1$ and $\sV$. All jet spaces we will use
are spaces of jets of sections of bundles over $M$. All vertical tangent vectors are vertical with respect to
the projection onto $M$. It means that $\sJ^1\sV E$ is the space of jets of sections of the bundle $\zz\circ
\zn_E$, therefore the projection onto $\sV E$ will be denoted by $\sj^1(\zz\circ \zn_E)$. Similarly, $\sV\sJ^1
E$ is the space of vectors tangent to $\sJ^1 E$ and vertical with respect to the projection onto $M$.  The
projection to $\sJ^1 E$ will be denoted by $\zn_{\sJ^1 E}$. Both spaces, $\sV\sJ^1 E$ and $\sJ^1\sV E$,  play
a very important role in the Lagrangian formulation of the classical field theory. We will discuss the structure of
these spaces in the next section. Now, we would like to point out that the coordinate system $(q^i, y^a)$ on
$V\subset E$ allows us to construct  coordinate systems:
$$(q^i, y^a, y^b{}_j, \delta y^c, \delta y^d{}_k)\qquad\text{on}\qquad (\sj^1\zz\circ\zn_{\sJ^1 E})^{-1}(V)\subset \sV\sJ^1 E$$
and
$$(q^i, y^a, \delta y^b, y^c{}_j, \delta y^d{}_k)\qquad\text{on}\qquad (\zn_E\circ \sj^1(\zz\circ\zn_E))^{-1}(V)\subset \sJ^1 \sV E.$$

\section{Lagrangian formulation}\label{seclagr}

In the first-order field theory, a Lagrangian is a map from the space of first jets of sections of the
bundle $\zz$ to scalar densities on $M$, covering the identity on $M$. Since we assumed that $M$ is oriented,
we can identify densities with $m$-covectors. A Lagrangian $L$ is therefore a map
$$L:\sJ^1 E\longrightarrow\Omega^m$$
covering the identity on $M$. The space $\sJ^1 E$ is often called {\it the space of infinitesimal
configurations} for the first-order field theory. Let us recall that in statics and other variational theories
all information about the system is contained in a constitutive set which is a subset of the cotangent bundle
of the configuration space. In mechanics we developed another description of a constitutive set, using so
called convenient representations of covectors. The convenient representation approach led to differential
equations for phase trajectories of a system. In field theory we can adopt the same scheme. In the infinitesimal
approach the space of infinitesimal configurations is a manifold and the role of the internal energy is played by
the Lagrangian. Since variations are vectors vertical with respect to the projection onto $M$, the constitutive
set is given as an image of the vertical differential of the Lagrangian. Taking into account that the
Lagrangian has values in $\Omega^{m}$, we get that the constitutive set is a subset of $\sV^\ast\sJ^1
E\otimes_{\sJ^1 E}\Omega^m$. In mechanics, studying convenient representations led us to the concept of
momentum and external force. Let us do the same for the field theory.

\subsection{The phase space}\label{secphase} The first step is recognizing the phase space for the first-order
field theory on the fibration $\zz$. We can use the calculus of variations as the guide-line for the problem.
Let $D$ be a compact region in $M$ with the smooth boundary $\partial D$ such that it is contained in a domain of
coordinates $U$. We will do now the standard calculations in coordinates that leads to the Euler-Lagrange
equations, but we will not assume that the variations vanish on the boundary. In the following, we will denote
by $(\zs^a)$ the functions defining local section $\sigma$, and by $(\zs^a, \zd\zs^a)$ the functions defining its
variation, i.e. a vertical vector field on $E$ along the section $\zs$. The action functional $S$ evaluated on
$\sigma$ gives
$$S[\sigma]=\int_D\ell(q^i, \zs^a, \frac{\partial\zs^b}{\partial q^j})\;\eta,$$
and the variation of $S$ evaluated on the variation of $\zs$ gives
\begin{equation}\label{lg1}
\langle\zd S,\zd\zs\rangle= \\
\int_D\left(\frac{\partial\ell}{\partial y^a}\,\zd\zs^a + \frac{\partial\ell}{\partial
y^b{}_j}\frac{\partial\zd\zs^b}{\partial q^j}\right)\eta.
\end{equation}
Using the Stokes theorem, we obtain
\begin{equation}\label{lg2}
\langle\zd S,\zd\zs\rangle= \int_D\left(\frac{\partial\ell}{\partial y^a}\zd\zs^a-\frac{\partial}{\partial
q^j}\frac{\partial\ell}{\partial y^a{}_j}\right)
\zd\zs^a\; \eta+\\
\int_{\partial D}\frac{\partial\ell}{\partial y^b{}_i}\,\zd\zs^b\;\eta_i.
\end{equation}
The term integrated over $D$ gives the definition of the external forces that are now interpreted as sources of a
field, while the term integrated over the boundary $\partial D$ gives the definition of the momenta. An object
that can be integrated over the boundary of the region $D$ is a $(m-1)$-form, therefore the results of the
evaluation of momenta over the variations should lie in the bundle of $(m-1)$-covectors on the base manifold
$M$. We point out, however, that the momenta should be evaluated on variations rather than on infinitesimal
configurations. It is a special case of autonomous mechanics when variations and infinitesimal configurations
are represented by the same geometrical object (tangent vector), therefore  we can evaluate momenta on
velocities. We cannot evaluate the momenta on first jets, at least not in the usual sense. The phase space for the
first-order field theory on the bundle $\zz$ is therefore the space
\begin{equation} \mathcal{P}=\sV^\ast E\otimes_E\Omega^{m-1}.\end{equation}
It is a vector bundle over $E$. We will denote the corresponding fibration by
$$\pi: \mathcal{P}\rightarrow E.$$
Using a base of sections of the bundle $\rho_E$ and a base $(\eta_i)$ in $\Omega^{m-1}$, we can construct
local linear coordinates on $\pi^{-1}(V)$,
\begin{equation}(q^i, y^a, p^j{}_b),\end{equation}
such that a section of the bundle $\pi$ is represented as
$$p^j{}_b(q,y)\xd y^b\otimes \eta_j.$$

The role of external forces is played by the source of a field. The source is represented by a section of the
bundle $V^\ast E\otimes_E\Omega^m\rightarrow M$.

All the above calculations can be done in a coordinate-free form. We postpone it to section \ref{phase1}.

\subsection{The structure of iterated bundles.}
Passing from (\ref{lg1}) to (\ref{lg2}) we have used implicitly a canonical mapping
$$\kappa:\sV\sJ^1 E\longrightarrow \sJ^1\sV E$$
which is analogous to the canonical flip $\kappa_M:\sT\sT M\rightarrow\sT\sT M$. {The role of $\kappa$ in the
field theory is similar to that of $\kappa_M$ in mechanics, where the map $\kappa_M$ is used to obtain
a convenient representation of a variation of an infinitesimal configuration. The idea of using convenient
representations comes from mechanics formulated for finite time intervals. In the field theory, instead of
the finite time intervals, we have the bounded domains of integration $D\subset M$. A configuration is then a section
of $\zz$ restricted to $D$ and its virtual displacement is represented by an equivalence class of curves in
the space of sections. Curves in the space of sections come from vertical homotopies, i.e. maps
$$\chi: U\times I\longrightarrow E$$
such that $D$ is contained in an open set $U$ and $I$ is a neighborhood of $0$ in $\R$. The verticality means
that for any $t$ we have $\zz(\chi(x,t))=x$. Fixing $x$ we obtain a vertical curve in $E_x$, while fixing $t$
we obtain a local section of $\zz$. Restricting the domain of $\chi$ to $D\subset U$, we obtain a curve in
configurations. The curves are classified, as usual, with a use of some functions on configurations of the
type of an action functional. We observe that equivalence classes are conveniently represented by vertical
vector fields $\delta\sigma$ along a sections $\sigma$ over $D$. In the infinitesimal approach, a configuration
is the first jet $\sj^1\sigma(x)$ of a section and its variation is represented by a vector $\delta\sj^1\sigma$
tangent to the space of first jets and vertical with respect to the projection on $M$. From the convenient
representation we obtain the first jet of a vertical vector field along a section, i.e. an element of $\sJ^1\sV
E$.}

Let us be more precise and define the map $\kappa$ using representatives of elements of $\sJ^1\sV E$ and
$\sV\sJ^1 E$. For a vertical homotopy
$$\chi: U\times I\longrightarrow E$$
of local sections of $\zz$, and for a point $x_0\in M$, we can create two objects.
Taking the first jet of the section $x\mapsto\chi(x,t)$ at $x_0$, we obtain the curve
$$t\mapsto \sj^1\chi(x_0,t)$$
in $\sJ^1 E$ which is vertical with respect to the projection onto $M$. The vector
$$\st\sj^1\chi(x_0,0),$$
tangent to this curve at $t=0$, is an element of $\sV\sJ^1 E$. On the other hand, we can first take vectors
tangent to vertical curves $t\mapsto \chi(x,t)$ at $t=0$, obtaining a vertical vector field along the section
$x\mapsto \chi(x,0)$. The vector field
$$x\mapsto\st\chi(x,0)$$
is a section of the bundle $\zeta\circ\rho_E:\sV E\rightarrow M$. Taking the first jet of this section at the
point $x_0$, we obtain the element $\sj^1\st\chi(x_0,0)$ of $\sJ^1\sV E$.
\begin{definition} The map $\kappa: \sV\sJ^1E\rightarrow \sJ^1\sV E$ is uniquely determined by
\be\label{kappa}\kappa(\st\sj^1\chi(x_0,0))=\sj^1\st\chi(x_0,0).\ee
\end{definition}
\noindent The definition is correct as the both sides of (\ref{kappa}) are independent on the choice of the
representative $\chi$ for an element of $\sV\sJ^1E$.
\smallskip

Both spaces, $\sJ^1\sV E$ and $\sV\sJ^1 E$, are double bundles in the sense that they carry the structure of two
compatible fibrations. The space $\sJ^1\sV E$ is fibrated over $\sV E$ and the fibration is an affine
bundle modeled on the vector bundle $\sT^\ast M \otimes_{\sV E} \sV\sV E\rightarrow \sV E$. The
projection from $\sJ^1\sV E$ onto $\sV E$ is $\sj^1(\zeta\circ\zn_E)$, since the first jets are calculated with
respect to the projection on $M$. In the adapted coordinates we have
$$(q^i, y^a, \delta y^b, y^c{}_j, \delta y^d{}_k)\longmapsto (q^i, y^a, \delta y^b).$$
The space $\sJ^1\sV E$ is also fibred over $\sJ^1 E$ and the corresponding projection will be denoted
$\sJ^1\nu_E$. In the adapted coordinates the second projection reads as
$$(q^i, y^a, \delta y^b, y^c{}_j, \delta y^d{}_k)\longmapsto (q^i, y^a, y^c{}_j).$$
The compatibility of the two bundle structures means that the model vector bundle $\sT^\ast M\otimes_{\sV
E}\sV\sV E$ is in fact a double vector bundle. For the concept of a double affine bundle and its model
double vector bundle we refer to \cite{GRU}.  We have therefore two commutative diagrams of bundle
projections
\begin{equation}\begin{array}{ccc}
\xymatrix@C-5pt{
 & \sJ^1\sV E\ar[dl]_{\sJ^1\zn_E}\ar[dr]^{\sj^1(\zz\circ\zn_E)} & \\
\sJ^1 E\ar[dr]^{\sJ^1\zeta} & & \sV E\ar[dl]^{\zn_E} \\
 & E &
} & \qquad & \xymatrix@!C=2pc{
 & \sT^\ast M \otimes_{\sV E}\sV\sV E\ar[dl]_{\sV\zn_E}\ar[dr]^{\zn_{\sV E}} & \\
\sT^\ast M\otimes_E\sV E\ar[dr]^{\zn_E} & & \sV E\ar[dl]_{\zn_E} \\
 & E &
}\end{array}
\end{equation}

The space $\sV\sJ^1 E$ carries also the structure of a double bundle. One of the bundles is an affine  bundle
and the other is a vector bundle. The vector bundle fibration $$\nu_{\sJ^1 M}: \sV\sJ^1
E\longrightarrow \sJ^1 E$$ in the adapted coordinates reads as
$$(q^i, y^a, y^b{}_j, \delta y^c, \delta y^d{}_k)\longmapsto (q^i, y^a, y^b{}_j).$$
The second fibration is the affine bundle fibration
$$\sV\sJ^1\zeta: \sV\sJ^1 E\longrightarrow \sV E$$
which we obtain applying the vertical tangent functor $\sV$ to the projection $\sJ^1\zeta:\sJ^1 E\rightarrow
E$. The model vector bundle for this affine bundle is
$$\sV (\sV E\otimes_E\sT^\ast M)\rightarrow \sV E.$$
In the adapted coordinates the affine projection $\sV\sJ^1\zeta$ reads as
$$(q^i, y^a, y^b{}_j, \delta y^c, \delta y^d{}_k)\longmapsto (q^i, y^a, \delta y^b).$$
The compatibility condition of the two projections can be expressed as the assumption that the model space
$\sV (\sV E\otimes_E\sT^\ast M)$ is a double vector bundle. The structure of $\sV\sJ^1 E$ can be summarized in
the following two commutative diagrams
\begin{equation}
\begin{array}{ccc}
\xymatrix@C-5pt{
 & \sV \sJ^1E\ar[dl]_{\zn_{\sJ^1 E}}\ar[dr]^{\sV\sj^1\zz} & \\
\sJ^1 E\ar[dr]^{\sJ^1\zeta} & & \sV E\ar[dl]_{\zn_E} \\
 & E &
} &\qquad & \xymatrix@!C=2pc{
 & \sV (\sV E\otimes_E\sT^\ast M)\ar[dl]_{\zn_{\sV E}}\ar[dr]^{\sV\zn_E} & \\
\sV E\otimes_E \sT^\ast M\ar[dr]^{\zn_E} & & \sV E\ar[dl]_{\zn_E} \\
 & E &
}  \end{array}
\end{equation}
where the second diagram represents the model double vector bundle.

The map $\kappa$ defined earlier is a morphism of double bundles covering the identities on side bundles. On the level of
the model double vector bundles it corresponds to the canonical flip $\kappa_E:\sT\sT E\rightarrow\sT\sT E$
restricted to vertical vectors.
\medskip

{Another example of a space equipped with the double structure of a vector-affine bundle is the space
$\sJ^1\mathcal{P}$ of first jets of the bundle $\zeta\circ\pi$. As the bundle of jets it is fibred over
$\mathcal{P}$ and the fibration is an affine fibration modeled on a vector fibration
$\sV\mathcal{P}\otimes_{\mathcal{P}}\sT^\ast M\rightarrow \mathcal{P}$. The vector bundle structure on
$\sJ^1\mathcal{P}$ comes from the jet prolongation of the vector bundle projection $\pi:\mathcal{P}\rightarrow
E$. The double bundle structure of $\sJ^1 E$ can be summarized in the following two diagrams,
\begin{equation}
\begin{array}{ccc}
\xymatrix@!C=2pc{
 & \sJ^1\mathcal{P}\ar[dr]^{\sJ^1\pi}\ar[dl]_{\sj^1(\zeta\circ\pi)} & \\
\mathcal{P}\ar[dr]^{\pi} & & \sJ^1 E\ar[dl]_{\sj^1\zz} \\
 & E &
} &\qquad & \xymatrix@!C=2pc{
 & \sV\mathcal{P}\otimes_{\mathcal{P}}\sT^\ast M\ar[dl]_{\nu_{\mathcal{P}}}\ar[dr]^{\sV\pi} & \\
\mathcal{P}\ar[dr]^\pi & & \sV E\otimes_E\sT^\ast M\ar[dl]_{\nu_E} \\
 & E &
}  \end{array}
\end{equation}
where the second diagram represents the model double vector bundle. Elements of the model vector bundle
$\sV\mathcal{P}\otimes_{\mathcal{P}}\sT^\ast M$ can be added to elements of $\sJ^1\mathcal{P}$. In
$\sV\mathcal{P}\otimes_{\mathcal{P}}\sT^\ast M$ there is a subbundle of  vectors which are vertical with
respect to the projection on $E$, i.e. vectors $v$ such that $\sV\pi(v)=0$. Vectors vertical with respect to
the projection on $E$ have the first component tangent to the corresponding fibre of the bundle $\pi$. As usual,
vectors tangent to a fibre of a vector bundle can be identified with elements of the fibre itself. Therefore, if $\sV\pi(v)=0$,
then $v$ can be identified with an element of $\mathcal{P}\otimes_E\sT^\ast M=\sV^\ast
E\otimes_E\Omega^{m-1}\otimes_E\sT^\ast M$. Note that adding vectors vertical with respect to the projection
on $E$ does not change the right-hand side projection, i.e., if $p$ is a local section of $\zz\circ\pi$, then
$$\sJ^1\pi(\sj^1p(x)+v)=\sJ^1\pi(\sj^1p(x)).$$
In coordinates, if $\sj^1p(x)=(x^i, y^a, p^j{}_b, y^c{}_k, p^l{}_{ds})$ and
$v=(x^i, y^a, v^j{}_{bk})$ i.e. $v=v^j{}_{bk}\xd y^b\otimes\eta_j\otimes\xd x^k$, then
$$\sj^1p(x)+v=(x^i, y^a, p^j{}_b, y^c{}_k, p^l{}_{ds}+v^l{}_{ds}).$$
In section \ref{secphase} we have observed that sources of fields are represented by sections of the bundle\newline
$\sV^\ast E\otimes_E \Omega^m\rightarrow M$. Since
$$\sV^\ast E\otimes_E \Omega^m\subset\sV^\ast E\otimes_E\Omega^{m-1}\otimes_E\sT^\ast M$$
we see that elements of the total space of the bundle of sources can be added to first jets from
$\sJ^1\mathcal{P}$ without changing any of the projections. This operation will be needed in construction of
the field phase equations with sources in section \ref{secalpha}.}

\bigskip
{In the next section we will construct the main map of the Lagrangian formulation of the field theory that
maps covectors on the space of infinitesimal configurations to their convenient representations.} For that we
will need an evaluation between the space $\sJ^1\mathcal{P}$ of first jets of sections of the bundle
$$\zeta\circ\pi:\;\mathcal{P}\longrightarrow M$$
and the space $\sJ^1\sV E$ of first jets of vertical virtual displacements. More precisely, we will construct
an evaluation between the bundle
$$\sJ^1\pi:\; \sJ^1 \mathcal{P}\longrightarrow\sJ^1 E$$
and the bundle
$$\sJ^1\zn_E:\; \sJ^1\sV E\longrightarrow\sJ^1 E$$
with values in the pull-back of the bundle $\Omega^m\rightarrow M$ by $\zz\circ\sj^1\zz$.
\smallskip

Let $p: M\supset U\longrightarrow \mathcal{P}$ be a local section of the momentum bundle. We denote by
$\sigma$ the underlying section of the bundle $\zeta$, i.e. $\sigma: M\supset U\longrightarrow E$ is  such
that $p\circ\pi=\sigma$. Let also $\delta\sigma : M\supset U\longrightarrow \sV E$ be a vertical vector field
along the section $\sigma$. There is a natural evaluation between $\sV_e E$ and $\mathcal{P}_e=\sV^\ast_e
E\otimes \Omega^{m-1}_x$ with values in $\Omega^{m-1}_x$, therefore $\langle p,\delta\sigma\rangle$ is a
$(m-1)$-form defined on $U\subset M$. W can define the evaluation between $\sj^1 p(x_0)$ and
$\sj^1\delta\sigma(x_0)$ using the formula
$$\langle\!\langle\,\sj^1 p(x_0),\sj^1\delta\sigma(x_0)\,\rangle\!\rangle=
\xd \langle p,\delta\sigma\rangle(x_0),$$ so that the evaluation is a map
$$\langle\!\langle\cdot,\cdot\rangle\!\rangle: \sJ^1\mathcal{P}\times_{\sJ^1 E}\sJ^1\sV E\longrightarrow \Omega^m.$$
In coordinates, if $\zs$ is given by local functions $(\sigma^a)$, if the momentum is represented as
$$p^i{}_a\xd y^a\otimes\eta_i$$
and $\delta\sigma$ as
$$\delta\sigma^b\frac{\partial}{\partial y^b},$$
then
$$\langle p,\delta\sigma\rangle(x)=p^i{}_a(x)\delta\sigma^a(x)\eta_i$$
and
$$\xd \langle p,\delta\sigma\rangle(x_0)=\left(\frac{\partial p^j{}_b}{\partial q^j}(x_0)\delta\sigma^b(x_0)+p^j{}_b(x_0)
\frac{\partial \delta\sigma^b}{\partial q^j}(x_0)\right)\eta,$$ therefore the evaluation $\langle\!\langle,
\rangle\!\rangle$ in coordinates reads as
$$\langle\!\langle (q^i, y^a, p^j{}_b, y^c{}_k, p^l{}_{dm} ),
(q^i, y^a, \delta y^b, y^c{}_k, \delta y^d{}_l) \rangle\!\rangle = p^l{}_{dl} \delta y^d+ p^j{}_b\delta
y^b{}_j.$$
\bigskip

\subsection{The map $\za$ }\label{secalpha}
Let us now define the main geometrical object of the Lagrangian formulation of the first-order field theory.
\begin{definition} The relation
$$\za: \sJ^1 \mathcal{P}\longrightarrow \sV^\ast\sJ^1 E\otimes_{\sJ^1 E}\Omega^m,$$
given by the condition
$$\langle\!\langle u, \kappa(w)\rangle\!\rangle =\langle \za(u), w\rangle$$
for all $w$ having the same projection on $\sJ^1 E $ as $u$, will be called the {\it Lagrangian relation}.
\end{definition}
\noindent In this case the relation $\za$ is actually a mapping. In coordinates we have
$$\za:\; (q^i,y^a,p^j{}_b, y^c{}_k, p^l{}_{dm})\longmapsto (q^i,y^a, y^c{}_k, \sum_lp^l{}_{dl}, p^j{}_b).$$
{The map $\za$ is a field-theoretical analog of the Tulczyjew $\za_M$ in mechanics. It relates covectors on
the space of infinitesimal configurations which are elements of $\sV^\ast\sJ^1
E\otimes_{\sJ^1 E}\Omega^m$, to their convenient representations. In the simplest case, when sources of the
field are equal to $0$, a convenient representation of a covector is the first jet of a section of the momentum
bundle.} If there are no constraints for infinitesimal configurations of the system described by the
Lagrangian $L$, the constitutive set of the system is given as an image of $\sJ^1 E$ by the vertical
differential $\xd L$ (understood as a map from $\sJ^1 E$ to $\sV^\ast \sJ^1 E\otimes_{\sJ^E}\Omega^{m}$). {
Using the map $\alpha$ we can obtain a convenient representation of the constitutive set as a differential
inclusion which we understand as a condition for sections of the momentum bundle. This differential inclusion
we will call the {\it phase dynamics} of the field.}

\begin{definition} The \emph{phase dynamics} of the field, when sources are equal to $0$,
is the subset $\mathcal{D}$ of $\sJ^1\mathcal{P}$ given by
$$\mathcal{D}=\alpha^{-1}(\xd L(\sJ^1 E)).$$
\end{definition}

\noindent {The phase dynamics is also called the {\it Lagrangian field equations.} Let us note that obtaining
the Lagrangian field equations from a Lagrangian is in our theory very simple. We do not require any regularity of
the Lagrangian.}
\begin{definition} We say that a section $p: M\rightarrow\mathcal{P}$ is a solution of
the Lagrange field equations if
$$\sj^1 p(x)\in \mathcal{D}.$$
\end{definition}
\noindent In coordinates it means that
$$\sum_j\frac{\partial p^j{}_b}{\partial q^j}=\frac{\partial L}{\partial y^b},\qquad p^j{}_b=\frac{\partial L}{\partial y^b{}_j}.$$
The equations, known as the Euler-Lagrange equations for field theory, are consequences of the Lagrange field
equations.
\smallskip

The Lagrangian side of the Tulczyjew triple for the first-order classical field theory can be presented in the
following diagram:
\begin{equation}
\xymatrix@!C=2pc{
 & \sJ^1\mathcal{P}\ar[rrr]^{\alpha}\ar[rd]^{\sJ^1 \pi}\ar[ldd]_/-10pt/{\sj^1(\pi\circ\zz)} & & &
 \sV^\ast\sJ^1 E\otimes_{\sJ^1 E}\Omega^m\ar[rd]^{\rho_{\sJ^1 E}}\ar[ldd]_/-10pt/{\xi} &  \\
 & & \sJ^1 E\ar[ldd]_/-10pt/{\sj^1\zz}\ar[rrr]^/-15pt/{id} & & & \sJ^1 E\ar[ldd]_/-10pt/{\sj^1\zz} \\
\mathcal{P}\ar[rd]^{\pi}\ar[rrr]^/+15pt/{id} & & & \mathcal{P}\ar[rd]^{\pi} & & \\
 & E\ar[rrr]^{id} & & & E &
}
\end{equation}
There is one projection in the above diagram that needs explanation. It is the projection
$$\zx: \sV^\ast \sJ^1 E\otimes_{\sJ^1 E}\Omega^m\longrightarrow\mathcal{P}.$$
Let us fix a point $v\in\sJ^1E$ and denote $e=\sj^1\zz(v)$ and $x=\zz(e)$. In the vector space $\sV_v\sJ^1E$
there is a subspace of those vectors which are tangent to fibres of the projection $\sj^1\zz:\sJ^1E\rightarrow
E$. Since each such fibre is an affine subspace of $\sT^\ast_x M\otimes\sT_e E$, the space of vectors tangent
to the fibre is isomorphic to its model vector space which is $\sT^\ast_x M\otimes\sV_e E$. An element of
$\sV^\ast_v \sJ^1 E\otimes\Omega^m_x$, treated as a linear function on $\sV_v \sJ^1 E$ with values in
$\Omega^m_x$, can be restricted to the subspace of vectors tangent to the fibres. The restriction is an
element of
$$\sT_x M\otimes\sV^\ast_e E\otimes\Omega^m_x\approx \sV^\ast_e E\otimes\Omega^{m-1}_x=\mathcal{P}_e.$$
Summarizing, the projection $\zx$ is a restriction of a covector to the subspace of vectors tangent to fibres
of a certain projection. It provides the Legendre map, defined by the Lagrangian, from the space of infinitesimal
configurations to the space of momenta,
$$\lambda: \sJ^1 E\longrightarrow\mathcal{P},\qquad\lambda(v)=\zx(\xd L(v)).$$

{Sources of the field can also be included in the picture. Recall that sources are sections of the bundle
$\sV^\ast E\otimes_E\Omega^m\rightarrow M$ and that an element of the total space of this bundle can be added
to the elements of $\sJ^1\mathcal{P}$. We define the extended map $\widetilde\alpha: \sV^\ast
E\otimes_E\Omega^m\times_{\mathcal{P}}\sJ^1\mathcal{P} \longrightarrow \sV^\ast\sJ^1E\otimes_{\sJ^1
E}\Omega^m$ by the formula
$$\widetilde\alpha(\rho, \sj^1p)=\alpha(\rho+\sj^1 p).$$

\begin{definition} The \emph{phase field dynamics with sources}
is the subset $\mathcal{\widetilde D}$ of $\sV^\ast\sJ^1E\otimes_{\sJ^1
E}\Omega^m\times_\mathcal{P}\sJ^1\mathcal{P}$ given by
$$\mathcal{\widetilde D}=\widetilde\alpha^{-1}(\xd L(\sJ^1 E)).$$
\end{definition}

\begin{definition} We say that a pair of sections $p: M\rightarrow\mathcal{P}$  and $\rho: M\rightarrow\sV^\ast E\otimes_{E}\Omega^m$ is a solution of phase field dynamics with sources if
$$\rho(x)+\sj^1 p(x)\in \mathcal{\widetilde D}.$$
\end{definition}
\noindent In coordinates it means that
$$\sum_j\frac{\partial p^j{}_b}{\partial q^j}+\rho_b=\frac{\partial L}{\partial y^b},\qquad p^j{}_b=\frac{\partial L}{\partial y^b{}_j}.$$}

\subsection{Phase space: geometrical version}\label{phase1}

In section \ref{secphase}, using a coordinate calculation, we have split the differential of a Lagrangian into
two parts: the Euler-Lagrange part and the total differential part. We have used the formula (\ref{lg2}) to
determine the phase space for the problem. Now we can do the same intrinsically, i.e. without using any specific
choice of coordinates (see \cite{Tu2}).

Let us fix a point $x_0$ in $M$. Any element $v$ of $\sV\sJ^1 E$ over $x_0$ can be represented as
$v=\kappa(\sj^1\delta\sigma(x_0))$ for some section $\delta\sigma$ of the bundle $\zz\circ\nu_E$. The same can
be done for any element of $\sV\sJ^l E$ with the use of the isomorphism
$$\kappa_{(l,1)}:\sJ^l\sV E\longrightarrow \sV\sJ^l E$$
which is defined analogously to $\kappa$, using representatives of elements of $\sV\sJ^l E$ and $\sJ^l\sV E$.
It is convenient to use the notation $\delta\sigma^l(x_0)=\kappa_{(l,1)}(\sj^l\delta\sigma(x_0))$. Given a
covector $\varphi\in\sT^\ast_{x_0}M$, we can choose a local function $f$ on $M$ such that $f(x_0)=0$ and $\xd
f(x_0)=\varphi$. Now we define
\begin{equation}
F(\varphi,v)=\kappa(\sj^1\zx(x_0)),\quad\text{where}\quad \zx(x)=f(x)\delta\sigma(x). \label{f1}
\end{equation}
It is clear that the value of $F$ depends only on the covector $\varphi$ and the vector $v$, and not on the
representatives. We have then defined the map
$$F: \sT^\ast_{x_0}M\times \sV_{\sj^1\sigma(x_0)}\sJ^1 E\longrightarrow \sV_{\sj^1\sigma(x_0)}\sJ^1 E$$
which is bilinear. In coordinates, if
$$\varphi=\varphi_i\xd q^i,\qquad v=\delta y^a\frac{\partial}{\partial y^a}+\delta y^a_j\frac{\partial}{\partial y^a_j},$$
then
$$F(\varphi, v)=\varphi_j\delta y^a\frac{\partial}{\partial y^a_j}.$$
We can see that the projection of $F(\varphi, v)$ on $\sV E$ is zero, i.e. $F(\varphi, v)$ is vertical with
respect to $\sj^1\zz$.

For a one-form $\zm$ on $\sJ^1 E$ with values in $\Omega^m$, we define a one-form $i_F\zm$ on $\sJ^1 E$ with
values in $\Omega^{m-1}$ by the formula
$$\langle\,\zm,\,F(\varphi,v)\,\rangle=\varphi\wedge\langle\, i_F\zm, v\rangle.$$
In coordinates, if
$$\zm=((\zm^0)_a\xd y^a+(\zm^1{})_a^i\xd y^a_i)\otimes\eta,$$
then
$$i_F\zm=(\zm^1{})_a^i\xd y^a\otimes\eta_i.$$
The form $i_F\zm$ is vertical, i.e. it vanishes on vectors vertical with respect to the projection $\sj^1\zz$.

There is an operation of total differential $\xd_M$ defined for forms on jet bundles with values in
$\Omega^k$. For example, if $\zm$ is a one-form on $\sJ^lE$ with values in $\Omega^k$, then $\xd_M\zm$ is a
one-form on $\sJ^{l+1}E$ with values in $\Omega^{k+1}$ given by the formula
$$\langle\,\xd_M\zm(\sj^{l+1}\sigma(x_0)),\, \delta\sigma^{l+1}(x_0)\,\rangle=
\xd\left(\langle\,\zm\circ\sj^l\sigma,\, \delta\sigma^l \rangle\right)(x_0).$$ Applying $\xd_M$ to $i_F\zm$,
for $\zm$ being a one-form on $\sJ^1 E$ with values in $\Omega^m$, we get that $\xd_Mi_F\zm$ is a one-form on
$\sJ^2E$ with values in $\Omega^m$. Let
$$E(\zm)=(\tau^2_1)^\ast\zm-\xd_Mi_F\zm\qquad\text{and}\qquad P(\zm)=i_F\zm,$$
where $\tau^2_1$ is the canonical projection
$$\tau^2_1: \sJ^2 E\longrightarrow \sJ^1E.$$
We have
$$(\tau^2_1)^\ast\zm=E(\zm)+\xd_M P(\zm)$$ and
both forms $E(\zm)$ and $P(\zm)$ are vertical with respect to the projection on $E$. In particular, if
$v\in\sV\sJ^2 E$ and $\sT\tau^2_E(v)=0$, then
$$\langle\, E(\zm),\, v\,\rangle=0.$$
Indeed, let us take a representative $\delta\sigma$ such that $v=\delta\sigma^2(x_0)$. Since $v$ is vertical,
$\delta\sigma(x_0)=0$, so
\begin{multline}\label{f2}
\langle\,E(\zm),\, v\rangle =
\langle\,\zm, \sT\tau^2_1(v) \rangle- \langle \xd_Mi_F,\, v\,\rangle=\\
\langle\,\zm, \delta\sigma^1(x_0)\,\rangle- \langle \xd_Mi_F,\, \delta\sigma^2(x_0)\,\rangle= \langle\,\zm,
\delta\sigma^1(x_0)\,\rangle- \xd\langle (i_F\zm)\circ\sj^1\sigma,\, \delta\sigma^1\,\rangle(x_0).
\end{multline}
The vector $\delta\sigma^1(x_0)$ is also vertical with respect to the projection on $E$, therefore we can find
a function $f$ vanishing at $x_0$ and a vector $u=\delta\omega^1(x_0)$ such that
$$\delta\sigma^1(x_0)=F(\xd f(x_0), u),$$
i.e. we can write (for the first part of formula (\ref{f2}))
$$\langle\,\zm, \delta\sigma^1(x_0)\,\rangle=
\langle\,\zm, F(\xd f(x_0), \delta\omega^1(x_0)\,\rangle= \xd
f(x_0)\wedge\langle\,i_F\zm(\sj^1\sigma(x_0)),\,\delta\omega^1(x_0)\rangle.$$ Using the fact that $f(x_0)=0$,
we can write that
$$\xd f(x_0)\wedge\langle\,i_F\zm(\sj^1\sigma(x_0)),\,\delta\omega^1(x_0)\rangle=
\xd\left( f\langle\,i_F\zm\circ\sj^1\sigma,\,\delta\omega^1\rangle\right)(x_0).
$$
Now, let us concentrate on the second part of formula (\ref{f2}). Since the form $i_F\zm$ is vertical, the
value of
$$\langle i_F\zm(\sj^1\sigma(x)), \delta\sigma^1(x)\rangle$$
depends only on the jet $\sj^1\sigma(x)$ of the base section and on $\delta\sigma(x)$. The value of the
differential
\begin{equation}\xd_M\left(\langle i_F\zm\circ\sj^1\sigma, \delta\sigma^1\rangle\right)(x_0)\label{f3}\end{equation}
depends therefore on the second jet $\sj^2\sigma(x_0)$ and the first jet $\delta\sigma^1(x_0)$. This means that
in  formula (\ref{f3}) we can substitute the section $\delta\sigma$ by the section $x\mapsto
\zx(x)=f(x)\delta\omega(x)$ that covers the same section $\sigma$ and has the same first jet at $x_0$. We can
now continue the calculation started in (\ref{f2}):
\begin{multline}\langle\,\zm(\sj^1\sigma(x_0)), \delta\sigma^1(x_0)\,\rangle- \xd\left(\langle i_F\zm\circ\sj^1\sigma,\, \delta\sigma^1\,\rangle\right)(x_0)= \\
\xd\left( f\langle\,i_F\zm\circ\sj^1\sigma,\,\delta\omega^1\rangle -  \langle i_F\zm\circ\sj^1\sigma,\,
\kappa(\sj^1\zx)\,\rangle\right)(x_0)= \xd\langle\,i_F\zm\circ\sj^1\sigma,\,f\delta\omega^1 -
\,\kappa(\sj^1\zx)\rangle(x_0).
\end{multline}
Let us now observe that $f(x)\delta\omega^1(x) - \,\kappa(\sj^1\zx(x))$ is a vertical vector for any $x$,
because $f(x)\delta\omega^1(x)$ projects onto $f(x)\delta\omega(x)$ and $\kappa(\sj^1\zx(x))$  projects on
$\delta\omega(x)$. Using verticality of $i_F\zm$, we see that
$$\langle\,i_F\zm\circ\sj^1\sigma,\,f\delta\omega^1 - \,\kappa(\sj^1\zx)\rangle$$
equals $0$ on the whole neighborhood of $x_0$, its differential is therefore equal to zero.\medskip

We have shown that, for any one-form $\zm$ on $\sJ^1 E$, the form $E(\zm)$ is vertical. The form $P(\zm)$ is
also vertical by definition. For $\zm=\xd L$ we can therefore define two maps:
$$\mathcal{E}(\xd L):\sJ^2E\longrightarrow \sV^\ast E\otimes_E\Omega^m,\qquad
\langle\mathcal{E}(L)(\sj^2\sigma(x)), \delta\sigma(x)\rangle=\langle\, E(\xd
L)(\sj^2\sigma(x)),\,\delta\sigma^2(x)\rangle,$$ and
$$\mathcal{P}(\xd L):\sJ^1E\longrightarrow \sV^\ast E\otimes_E\Omega^{m-1},\qquad
\langle \mathcal{P}(L)(\sj^1\sigma(x)),\, \delta\sigma(x)\rangle= \langle E(\xd
L)(\sj^1\sigma(x)),\,\delta\sigma^1(x)\rangle,$$ such that formula (\ref{lg2}) assumes the form
$$\langle\zd S,\zd\zs\rangle=\int_D \langle\mathcal{E}(\xd L)(\sj^2\sigma(x)), \delta\sigma(x)\rangle+
\int_{\partial D} \langle \mathcal{P}(\xd L)(\sj^1\sigma(x)),\, \delta\sigma(x)\rangle.$$

\subsection{Simple example: electrostatics}
Let us write the Lagrangian side of the Tulczyjew triple for Electrostatics. In this particular example we
have included also sources of the field. Since it is a nonrelativistic theory, we have decided to use a
three-dimensional affine space $A$ as our playground. For the affine space tangent and cotangent bundles are
trivial, so it is possible to present mathematical objects that appear in the theory in a simple way. The
model vector space for the affine space $A$ will be denoted with $V$. It is a three-dimensional vector space
equipped with a symmetric non-degenerate positive definite bilinear two-form $g$ representing the metric. For
the purpose of this example we do not assume that $A$ is oriented. For the affine space $A$ we have trivial
tangent and cotangent bundles,
$$\tau_A:A\times V\longrightarrow A, \qquad \pi_A:A\times V^\ast\longrightarrow A.$$
Moreover, because of the presence of the metric, there is a canonical isomorphism
$$\tilde g:V\longrightarrow V^\ast, \qquad v\longmapsto g(v,\cdot),$$
and a canonical scalar density $\mathfrak{g}\in\bigwedge^3_{o}V^\ast$. The symbol $\bigwedge^k_{o}V^\ast$
denotes the space of odd $k$-forms on $V$.

The potential of the electrostatics is a scalar field, therefore we take $E=A\times\R$ and $M=A$. The main
bundle of the theory is the trivial bundle
$$pr_A: A\times\R\longrightarrow A.$$
The space of the first jets of sections of the above bundle can be identified with
$$\sJ^1E=A\times\R\times V^\ast.$$
The first jet of a section $\varphi$ at a point $x$ is just $(x, \varphi(x), \xd\varphi(x))$. Since in our case
$\Omega^m(M)=A\times\bigwedge^3_{o}V^\ast$, we get that the Lagrangian is a map
$$L:A\times\R\times V^\ast\longrightarrow A\times\bigwedge^3{}_{o}V^\ast$$
covering the identity on $A$ which reads as
$$L(x,r,\mu)=\frac12\langle\mu, \tilde g^{-1}(\mu)\rangle\mathfrak{g}.$$
If $\varphi$ is a section of $pr_A$ (i.e. a function on $A$), then
$$L(x,\varphi(x),\xd\varphi(x))=\frac12\langle\xd\varphi(x), \tilde g^{-1}(\xd\varphi(x))\rangle\mathfrak{g}.$$

Let us look at the other spaces involved in the theory, namely
$$\sV^\ast\sJ^1E\otimes\Omega^m(M)=A\times\R\times V^\ast\times(\R\otimes \bigwedge^3{}_{o}V^\ast)\times
(V\otimes\bigwedge^3{}_{o}V^\ast)\simeq A\times\R\times V^\ast\times\bigwedge^3{}_{o}V^\ast\times
\bigwedge^2{}_{o}V^\ast.$$ An element of the above space will be denoted by
$(x,r,\mu,\eta,\vartheta).$
There are two projections:
$$\nu:\, A\times\R\times V^\ast\times\bigwedge^3{}_{o}V^\ast\times\bigwedge^2{}_{o}V^\ast\longrightarrow\ A\times\R\times V^\ast$$
and
$$\xi:\,   A\times\R\times V^\ast\times\bigwedge^3{}_{o}V^\ast\times\bigwedge^2{}_{o}V^\ast\longrightarrow\
 A\times\R\times\bigwedge^2{}_{o}V^\ast.$$
The graph of the differential of the Lagrangian is a subset defined as
$$\xd L(A\times\R\times V^\ast)=\{(x,r,\mu, 0, \,\tilde g^{-1}(\mu)\lrcorner \mathfrak{g}\,),\;\ x\in A,\, r\in\R,\, \mu\in V^\ast\}.$$
We see from the above that the phase space of the theory is
$$\mathcal{P}=A\times \R\times \bigwedge^2{}_{o}V^\ast.$$
We are looking for an equation for a section of the phase bundle over $A$, i.e. for a map
\be A\ni x\longmapsto (x, \varphi(x), E(x)).\label{elec1}\ee
The Legendre map that associates the phase to a configuration is
$$\lambda:A\times \R\times V^\ast\longrightarrow A\times\R\times\bigwedge^2{}_{o}V^\ast, \qquad
\lambda(x,r,\mu)=(x,r,\,\tilde g^{-1}(\mu)\lrcorner \mathfrak{g}\,).$$ Elements of $\bigwedge^2_{o}V^\ast$,
i.e odd two-forms, can also be interpreted as vector densities. In the case of electrostatics, the two-form or
the vector density that we obtain here, integrated over a surface in $A$, gives the flux of the electrostatic
field through the surface.

The space $\sJ^1 \mathcal{P}$ is in our case
$$\sJ^1\mathcal{P}=A\times\R\times \bigwedge^2{}_{o}V^\ast\times V^\ast\times(V^\ast\otimes \bigwedge^2{}_{o}V^\ast).$$
An element of the above space will be denoted by
$(x,r,p,\mu,\nu).$
Using the same symbols $(x,r,\mu)$ does not lead to any confusion, because the objects denoted by those
symbols are conserved by every map we use. The analog of the space of external forces for electrostatic field
is
$$\sV^\ast E\otimes\Omega^m(M)=A\times\R\times \bigwedge^3{}_{o}V^\ast, $$
whose elements are $(x,r,\rho)$. The extended map $\widetilde\alpha$ reads as
$$A\times\R\times \bigwedge^3{}_{o}V^\ast\times \bigwedge^2{}_{o}V^\ast\times V^\ast\times (V^\ast\otimes\bigwedge^2{}_{o}V^\ast)
\longrightarrow A\times\R\times V^\ast\times \bigwedge^3{}_{o}V^\ast\times \bigwedge^2{}_{o}V^\ast,$$
$$(x,r,\rho,p,\mu,\nu)\longmapsto (x,r,\mu, \text{Alt}(\nu)-\rho, p).$$
It means that the inverse image of $\xd L(A\times\R\times V^\ast)$ by $\alpha$ is the set $D$ of all
$(x,r,\rho,p,\mu,\nu)$ such that
\begin{equation}\label{elect4}
p=\tilde g^{-1}(\mu)\lrcorner \mathfrak{g},\quad \text{Alt}(\nu)=\rho.
\end{equation}
The set $D$ represents an equation for the section (\ref{elec1}) of the phase bundle: if $r=\varphi(x)$ and
$p=E(x)$, we get
$$ \mu=\xd\varphi(x),\quad \nu=\sj^1 E(x),$$
and the equations are
\begin{eqnarray} & E(x)=\tilde g^{-1}(\xd\varphi(x))\lrcorner \mathfrak{g},\label{elec2} \\ & \xd E(x)=\rho.\label{elec3}\end{eqnarray}
It is easy to see that substituting (\ref{elec2}) to (\ref{elec3}) we obtain
$$ \xd(\tilde g^{-1}(\xd\varphi(x))\lrcorner \mathfrak{g})=\rho$$
which can be written as
\be (\Delta\varphi)\mathfrak{g}=\rho,\ee
i.e the Poisson equation for the potential of electrostatic field produced by the charge density $\rho$.

All the mathematical objects of the theory have clear physical meaning: the field itself is a potential for an
electrostatic field. The electrostatic field is a covector field rather than a vector field, however in the
presence of metric we can canonically translate one into the other. The phase is a vector density associated
to the field and used to calculate the flux of the electrostatic field through a surface.

\section{Hamiltonian formulation}\label{secham}

In this section we will construct the Hamiltonian side of the Tulczyjew triple. {The name "Hamiltonian" is
usually associated with the time evolution of the system. In our approach, the Hamiltonian side of the triple
gives just another way of generating phase dynamics. The interpretation of the Hamiltonian itself strongly
depends on the particular theory.

The space of infinitesimal configurations $\sJ^1E$ is an affine bundle over $E$, therefore it is necessary to
use the affine geometry and the notion of affine duality. In the following we will recall this notion and
construct an affine analog of the cotangent bundle. Since the affine bundle that appears in field theory, i.e.
$\sj^1\zz: \sJ^1E\rightarrow E$, has a reach internal structure, we decided to work first with simpler objects,
and then apply the results to $\sj^1\zz$.}

\subsection{The affine-dual bundle}

Let us first recall some facts from the geometry of affine spaces (for more details see e.g \cite{GGU1}). Let
$\tau:A\rightarrow N$ be an affine bundle modeled on a vector bundle $\nu: V\rightarrow N$. Let us also fix a
one-dimensional vector space $U$. We will use the symbol $\nu^\ast$ for the projection $V^\ast\otimes_N
U\rightarrow N$.

The vector space of all affine maps from a fibre $A_q$, $q\in N$, to $U$ will be denoted by $\Aff(A_q,U)$.
Every affine map has its linear part, therefore $\Aff(A_q,U)$ is fibrated over $\Lin(V_q, U)\simeq V^\ast_q
\otimes U$. Collecting the spaces $\Aff(A_q,U)$ point by point in $N$, we obtain a vector bundle
\begin{equation}\label{hm3.1}\tau^\dag: \Aff(A,U)\rightarrow N\end{equation}
and an affine bundle
\begin{equation}\label{hm3}
\theta: \Aff(A,U)\longrightarrow V^\ast\otimes_N U.
\end{equation}
The fibration $\theta$ is an affine bundle modeled on the trivial vector bundle
$$pr_1:V^\ast\otimes_N U\times U\longrightarrow V^\ast\otimes_N U.$$
In the case $U=\R$, the space $\Aff(A_q,\R)$ is called  the {\it affine dual} of $A_q$ and denoted
$A_q^\dag$.
\smallskip

It is always useful to write geometrical objects in coordinates. We will use a set of coordinates adapted to
the structure. Let $(x^i)$ denotes a local system of coordinates in $\mathcal{O}\subset N$. Choosing a local
basis $e=(e_\alpha)$ of sections of $V$, the dual basis $\epsilon=(\epsilon^\alpha)$ of sections of $V^\ast$,
a reference section $a_0:\mathcal{O}\rightarrow A$, and a non-zero vector $u\in U$, we can construct the
adapted system of coordinates,
$$(x^i, f^\alpha)\quad\text{in}\quad A\quad\text{such that}\quad f^\alpha(a)=\epsilon^\alpha(a-a_0(q)),\; q=\tau(a),$$
and the adapted system of coordinates
\begin{equation}\label{hm5}
(x^i, \varphi_\alpha, r)\quad\text{in}\quad \Aff(A,U)
\end{equation}
such that
$$\varphi(a)=(\varphi_\alpha\epsilon^\alpha(a-a_0(q))+r)u,\quad q=\tau(a)$$
for  $\varphi\in \Aff(A,U)$. In coordinates, the projection $\theta$ is expressed as the projection onto first
two sets of coordinates $(x^i, \varphi_\alpha)$.
\smallskip

In the family of all smooth maps from $\Aff(A,U)$ to $U$ we distinguish maps $\Psi:\Aff(A,U)\rightarrow U$
which are affine along fibres of $\theta$ and satisfy
\begin{equation}\label{hm1}
\Psi(\varphi+u)=\Psi(\varphi)-u.
\end{equation}
Property (\ref{hm1}) implies that in every fibre $\theta^{-1}(p)$ of the fibration $\theta$ there is exactly
one point $\varphi_p$ such that $\Psi(\varphi_p)=0$. It means that the set $\Psi^{-1}(0)$ is the graph of a
section $\Sigma_\Psi$ of the fibration $\theta$. On the other hand, having a section $\Sigma$ we can define
$$\Psi_\Sigma(\zf)=\Sigma(\theta(\zf))-\zf.$$
It is obvious that $\Psi_{\Sigma}$ satisfies condition (\ref{hm1}). Therefore we have a one-to-one
correspondence between smooth maps which are affine along fibres and satisfy condition (\ref{hm1}) on one
hand, and smooth sections of the bundle $\theta$ on the other.
\smallskip

The differentials of maps satisfying property (\ref{hm1}) are such covectors on $\Aff(A,U)$ with values in $U$
that, restricted to vectors tangent to fibres of $\theta$, they give $-id$. The submanifold of such covectors
will be denoted $K_{-id}$. It is a coisotropic submanifold of $\sT^\ast \Aff(A,U)\otimes_{\Aff(A,U)} U$ with
respect to the canonical symplectic structure on $\sT^\ast \Aff(A,U)\otimes_{\Aff(A,U)} U$ with values in $U$.
Using the local system of coordinates (\ref{hm5}), we can construct the adopted system of coordinates on
$\sT^\ast \Aff(A,U)\otimes_{\Aff(A,U)} U$:
$$(x^i, \varphi_\alpha, r, \sigma_i, f^\alpha,\rho).$$
In the above coordinates the submanifold $K_{-id}$ is given by the condition $\rho=-1$ and the canonical
symplectic structure on $\sT^\ast \Aff(A,U)\otimes_{\Aff(A,U)} U$ reads as
$$\omega_{\Aff(A,U)}=\xd \sigma_i\wedge \xd x^i+\xd f^\alpha\wedge\xd \varphi_\alpha+\xd\rho\wedge\xd r.$$
Having a coisotropic submanifold, we can perform a symplectic reduction. Leaves of characteristic foliation are
orbits of the cotangent lift of the natural action of $U$ on $\Aff(A,U)$. The reduced manifold will be denoted
by $\sP \Aff(A,U)$. In coordinates the reduction is the map
$$K_{-id}\ni (x^i, \varphi_\alpha, r, \sigma_i, f^\alpha, -1)\longmapsto
(x^i, \varphi_\alpha, \sigma_i, f^\alpha)\in \sP \Aff(A,U).$$ Elements of $\sP \Aff(A,U)$ can be interpreted
also as equivalence classes of sections of the bundle $\theta$ with respect to the following equivalence
relation. Since $\Aff(A,U)$ is fibrated over $V^\ast\otimes_N U$ and the fibration is modeled on trivial
fibration with the fibre being $U$, the difference $\Sigma_2-\Sigma_1$ of two sections is a map from
$V^\ast\otimes_N U$ to $U$. It is therefore clear what means that $\xd(\Sigma_2-\Sigma_1)(p)=0$ for some
$p\in V^\ast\otimes U$. We say that two pairs
$$(p_1, \Sigma_1) \text{ and }\;(p_2, \Sigma_2) \;\text{are equivalent if and only if}\;\;
p_1=p_2 \;\text{and}\; \xd(\Sigma_2-\Sigma_1)(p_1)=0.$$ The equivalence class of $(p,\Sigma)$ is sometimes
denoted by $\xd\Sigma(p)$ and called the {\it differential of the section $\Sigma$ at the point $p$}. The manifold
$\sP\Aff(A,U)$ is obviously fibrated over $V^\ast\otimes_N U$. The fibration is an affine bundle modeled on
$\sT^\ast(V^\ast\otimes_N U)\otimes_{\Aff(A,U)} U\rightarrow V^\ast\otimes_N U$.

The above construction of $\sP\Aff(A,U)$ is analogous to the construction of $\sP Z$ for an affine bundle
$Z\rightarrow M$, modeled on the trivial bundle $M\times \R\rightarrow M$ given in \cite{GGU1}. The only
difference is that $\R$ is replaced by the one-dimensional vector space $U$.
\smallskip

The affine bundle $\sP\Aff(A,U)$ is actually a double bundle. The second bundle structure is inherited from the
double vector bundle $\sT^\ast \Aff(A,U)\otimes_{\Aff(A,U)} U$. Let us first recall that the structure of the
double vector bundle $\sT^\ast \Aff(A,U)\otimes_{\Aff(A,U)} U$
\begin{equation}\label{hm8}
\xymatrix@!C=3pc{& \sT^\ast \Aff(A,U)\otimes_{\Aff(A,U)} U\ar[dl]_{\pi_{\Aff(A,U)}}\ar[dr]^{\varsigma} & \\
\Aff(A,U)\ar[dr]^{\tau^\dag} & & \Aff(A,U)^\ast\otimes_N U\ar[dl]_{(\tau^\dag)^\ast} \\
& N &  }\end{equation} is given by the two commuting Euler vector fields associated with the two vector bundle
structures: the canonical one,
$$\nabla_1=\sigma_i\frac{\partial}{\partial\sigma_i}+ f^\alpha\frac{\partial}{\partial f^\alpha}+ \rho\frac{\partial}{\partial\rho},$$
and the second one,
$$\nabla_2=\varphi_\alpha\frac{\partial}{\partial\varphi_\alpha}+ f^\alpha\frac{\partial}{\partial f^\alpha}+ r\frac{\partial}{\partial r}.$$
The second projection $\varsigma$ can be understood as follows. The covector
$\psi\in\sT^\ast_\varphi\Aff(A,U)\otimes U$ can be restricted to vectors tangent at $\varphi$ to the fibre of
the projection $\Aff(A,U)\rightarrow N$, i.e. to the space $\Aff(A_q, U)$. Any such fibre is a vector space,
therefore vectors tangent to the fibre can be identified with elements of the fibre itself. This leads to the
identification of the restriction of $\psi$ with an element of the dual to the fibre, i.e. an element of
$\Aff(A_q, U)^\ast\otimes U$. We denote by $(\tau^\dag)^\ast$ the projection from $\Aff(A, U)^\ast\otimes_N U$
to $N$.

The coisotropic submanifold $K_{-id}$ is an affine subbundle of the canonical bundle structure
$\pi_{\Aff(A,U)}$ and a vector subbundle (over a submanifold, see \cite{GR}) of the second  bundle structure
$\varsigma$. The image $\varsigma(K_{-id})$ consists of all elements $h$ of $\Aff(A_q, U)^\ast\otimes U$ that
satisfy property (\ref{hm1}):
$$h(\varphi+u)=h(\varphi)-u.$$
Elements of $\Aff(A_q, U)^\ast\otimes U$ satisfying property (\ref{hm1}) are in a one-to-one
correspondence with elements of $A_q$ itself. A natural identification can be established as follows. Every
$a\in A_q$ gives rise to a linear map on $\Aff(A_q,U)$ with values in $U$ by evaluation, i.e.
$$h_a: \varphi\longmapsto \varphi(a),$$
however $h_a$  does not satisfy (\ref{hm1}). Property (\ref{hm1}) is satisfied by $-h_a$, i.e. the map
$\varphi\mapsto -\varphi(a)$.

For dimensional reasons, every element of $\Aff(A_q, U)^\ast\otimes U$ that satisfies  property (\ref{hm1}) is
of the form $-h_a$ for some $a\in A_q$.
\smallskip

The Euler vector field $\nabla_2$ is tangent to $K_{-id}$  and projectable with respect to the symplectic
reduction. It gives rise to a vector bundle structure
$$\sP\Aff(A,U)\rightarrow A.$$
We have therefore the double bundle
\begin{equation}\label{hm6}
\xymatrix@!C=2pc{& \sP \Aff(A,U)\ar[dl]_{\sP\theta}\ar[dr]^{\sP\varsigma} & \\
V^\ast \otimes_N U\ar[dr]^{\nu^\ast} & & A\ar[dl]_{\tau} \\
& N &  }
\end{equation}
with the left projection being an affine bundle and the right projection being a vector bundle.
\bigskip

We can repeat the above constructions for $A=(\sj^1\zeta\circ\zeta)^{-1}(x)$, i.e. the fibre of the bundle
$\sJ^1 E\rightarrow M$ that is fibred over $N=E_x$, $V=(\zeta\circ\rho_E)^{-1}(x)$, and $U=\Omega^m_x$. For
simplicity, we will denote $(\sj^1\zeta\circ\zeta)^{-1}(x)$ with $(\sJ^1_xE)$. The space
$\Aff(\sJ^1_eE,\Omega^m_x)$ will be called the {\it affine dual} of $\sJ^1_e E$ and denoted by $\sJ^\dag_e E$.
Usually the affine dual of an affine space is the vector space off all affine functions on the affine space
with real values. Here we replace real numbers with top forms on $M$, but we keep the name. As a result of the
above general construction, we get the affine dual $\sJ^\dag_x E$ together with the affine fibration
$$\theta_x : \sJ^\dag_x E\longrightarrow\mathcal{P}_x$$
and the correspondence between sections of $\theta_x$ and maps which are affine along fibres of $\theta_x$
with values in $\Omega^m_x$. Collecting affine dual spaces $\sJ^\dag_x E$ point by point in $M$, we obtain a
vector bundle
$$\sj^\dag\zz: \sJ^\dag E\longrightarrow E$$
and an affine bundle
$$\theta: \sJ^\dag E\longrightarrow \mathcal{P}.$$
Sections of the fibration $\theta$ are in a one-to-one correspondence with maps from $\sJ^\dag E$ to
$\Omega^m$ covering the identity on $M$, affine along fibres of $\theta$, and satisfying property (\ref{hm1}).
We have also the manifold $\sP \sJ^\dag E$ fibred over $\mathcal{P}$ and equipped with the canonical family of
symplectic forms with values in $\Omega^m$, parameterized by points in $M$ and obtained by reduction from
$V^\ast\sJ^\dag E\otimes_{\sJ^\dag E} \Omega^m$.

{We then expect that the Hamiltonian description of the first-order field theory will be connected with the
fibration $\theta:\sJ^\dag E\rightarrow \mathcal{P}$ and the space $\sP\sJ^\dag E$. In particular, a
Hamiltonian is a section of $\theta$ and the differential of this Hamiltonian at $p\in\mathcal{P}$ is an
element of $\sP\sJ^\dag E$, therefore the image of the differential of a Hamiltonian is a subset of
$\sP\sJ^\dag E$. } Diagram (\ref{hm6}) for $\sJ^1 E$ takes the form
\begin{equation}\label{hm7}
\xymatrix@C-10pt{& \sP \sJ^\dag E\ar[dl]_{\sP\theta}\ar[dr]^{\sP\varsigma} & \\
\mathcal{P}\ar[dr]^{\pi} & & \sJ^1 E \ar[dl]_{\sj^1\zeta}\\
& E &  }
\end{equation}

Let us end this section by constructing coordinates in $\sP\sJ^\dag E$ adapted to the structure of the double
bundle. In section \ref{sec2} we have introduced coordinates $(q^i, y^a, y^b{}_j)$ in $\sJ^1E$ and in section
\ref{secphase} coordinates $(q^i, y^a, p^j{}_a)$ in $\mathcal{P}$. Since $\sJ^\dag E$ is fibred over
$\mathcal{P}$ and the fibration is an affine bundle modeled on trivial bundle
$\mathcal{P}\times_M\Omega^m\rightarrow \mathcal{P}$, it will be convenient to use coordinates $(q^i, y^a,
p^j{}_a, r)$ in $\sJ^\dag E$ such that $r$ is an affine coordinate along the fibres of $\theta$. If
$\varphi=(q^i, y^a, p^j{}_b, r)$ and $\sj^1\sigma(x)=(q^i, y^a, y^b{}_j)$, then
$$\varphi(\sj^1\sigma(x))=p^j{}_by^b{}_j+r.$$
In $\sV^\ast\sJ^\dag E\otimes_{\sJ^\dag E}\Omega^m$ we have therefore the adopted coordinate system $(q^i,
y^a, p^j{}_a, r, \xi_a, y^b{}_k, \rho)$ and the coisotropic submanifold $K_{-id}$ is given by the condition
$\rho=-1$. Diagram (\ref{hm8}) in the case $A=\sJ^1 E$ takes the form
\begin{equation}\label{hm9}
\xymatrix@!C=3pc
{& V^\ast\sJ^\dag E\otimes_{\sJ^\dag E}\Omega^m  \ar[dl]_{\pi_{\sJ^\dag E}}\ar[dr]^{\varsigma} & \\
\sJ^\dag E\ar[dr]^{\sJ^\dag\zeta} & & (\sJ^\dag E)^\ast\otimes_E\Omega^m \ar[dl]_{(\sJ^\dag\zeta)^\ast}\\
& E &  }
\end{equation}
that in coordinates reads as
\begin{equation}\label{hm11}
\xymatrix@!C=3pc
{& (q^i, y^a, p^j{}_a, r, \xi_a, y^b{}_k, \rho)  \ar[dl]_{\pi_{\sJ^\dag E}}\ar[dr]^{\varsigma} & \\
(q^i, y^a, p^j{}_a, r)\ar[dr]^{\sJ^\dag\zeta} & & (q^i, y^a, y^b{}_j, \rho) \ar[dl]_{(\sJ^\dag\zeta)^\ast}\\
& (q^i,y^a) &  }
\end{equation}
After the reduction, we obtain coordinates $(q^i, y^a, p^j{}_a,\xi_a, y^b{}_k)$ in $\sP\sJ^\dag E$. Diagram
(\ref{hm7}) in coordinates takes the form
\begin{equation}\label{hm10}
\xymatrix@!C=3pc
{& (q^i, y^a, p^j{}_a,\xi_a, y^b{}_k)\ar[dl]_{\sP\theta}\ar[dr]^{\sP\varsigma} & \\
(q^i, y^a, p^j{}_a)\ar[dr]^{\pi} & & (q^i, y^a, y^b{}_k)\ar[dl]_{\sj^1\zeta}\\
& (q^i, y^a) &  }
\end{equation}
The family of symplectic forms with values in $\Omega^m$ parameterized by points of $M$ is given in
coordinates as
\begin{equation}\label{hm12}
\omega_{\sP\sJ^\dag E}=(\xd\xi_a\wedge\xd y^a+\xd y^b{}_k\wedge\xd p^k{}_b)\otimes \eta.
\end{equation}

\subsection{The map $\beta$}
{In this section we will construct the map $\beta$ which will be used in deriving the phase dynamics of the
field from a Hamiltonian. It will be a field-theoretical version of $\beta_Q$ (see
(\ref{intro8})).} In mechanics there is a well-known formula relating Lagrangians to Hamiltonians, namely
\begin{equation}\label{be1}
    H(p)=\langle p,v\rangle-L(v).
\end{equation}
The origin of this formula lies in the procedure of composing symplectic relations \cite{ST,Tu8,B}. Recall
that a symplectic relation between symplectic manifolds $(P_1,\omega_1)$,  $(P_2,\omega_2)$ is a Lagrangian
submanifold in $(P_1\times P_2, \omega_1-\omega_2)$. If we deal with cotangent bundles, we can think of
generating objects for symplectic relations. For example, it is well known that there is a canonical
symplectomorphism $\mathcal{R}_F:\sT^\ast F\rightarrow\sT^\ast F^\ast$ for any vector bundle
$\tau:F\rightarrow M$. The graph of $\mathcal{R}_F$ is the Lagrangian submanifold in $(\sT^\ast
F\times\sT^\ast F^\ast, \omega_F-\omega_{F^\ast})$ generated by the evaluation of covectors and vectors
$$F\times_M F^\ast\ni(f,\varphi)\longmapsto \varphi(f)\in\R$$
in the following sense. The evaluation is a function defined on the submanifold $F\times_M F^\ast\subset
F\times F^\ast$, therefore it generates a Lagrangian submanifold in $\sT^\ast(F\times F^\ast)$. The space
$\sT^\ast(F\times F^\ast)$ can be naturally identified with $(\sT^\ast F\times\sT^\ast F^\ast,
\omega_F+\omega_{F^\ast})$. To get a Lagrangian submanifold with respect to the form
$\omega_F-\omega_{F^\ast}$, we apply the transformation
$$\sT^\ast F\times\sT^\ast F^\ast\ni(\zeta_1,\zeta_2)\mapsto(\zeta_1,-\zeta_2)\in\sT^\ast F\times\sT^\ast F$$
to the generated submanifold. In the adapted coordinates $(x^i, f^\alpha)$ in $F$,
$(x^i,f^\alpha,\sigma_i,\psi_\alpha)$ in $\sT^\ast F$, and $(x^i,\varphi_\alpha,\sigma_i,h^\alpha)$ in
$\sT^\ast F^\ast$, the isomorphism $\mathcal{R}_F$ reads
$$\mathcal{R}_F(x^i,f^\alpha,\sigma_i,\psi_\alpha)=(x^i,\psi_\alpha,\sigma_i,-f^\alpha).$$
If $L$ is a function on $F$, the Lagrangian submanifold $\xd L(F)\subset\sT^\ast F$ can be treated as a
symplectic relation between the cotangent bundle of the trivial one-point manifold  and $\sT^\ast F$. As a
symplectic relation the submanifold is generated by $-L$. The operation of composition of symplectic relations
does not in general lead to a symplectic relation, but since $\mathcal{R}$ is a diffeomorphism, we do not
encounter such problems here. The set $\mathcal{R}_F(\xd L(F))$ is a Lagrangian submanifold in $\sT^\ast
F^\ast$. Composing symplectic relations means adding generating objects, we have therefore a generating family
\begin{equation}\label{be2}
H: F^\ast\times_M F\longrightarrow \R,\qquad H(f,\varphi)=\varphi(f)-L(f),
\end{equation}
that in some cases can be reduced to a single function on $F^\ast$.

To clarify all sign problems let us mention that the family (\ref{be2}) is a generating object for a
symplectic relation between a single point and $\sT^\ast F^\ast$. For a generating object of the Lagrangian
submanifold $\mathcal{R}_F(\xd L(F))$ we have to take $-H$.

In mechanics the whole procedure is applied to $F=\sT Q$. Then, $-H$ is the generating object of a Lagrangian
submanifold $D_H$ in $\sT^\ast\sT^\ast Q$. The dynamics, i.e. a subset of $\sT\sT^\ast Q$, is obtained as the
inverse image of $D_H$ by the map $\beta_Q:\sT\sT^\ast Q\rightarrow \sT^\ast\sT^\ast Q$,
\begin{equation}\label{be3}\beta_Q=\alpha_Q\circ\mathcal{R}_{\sT Q}.\end{equation}
In this case $\beta_Q$ is also associated with $\omega_Q$
\begin{equation}\label{be4}\beta_Q(v)=\omega_Q(v, \cdot)\quad\text{for}\quad v\in\sT\sT^\ast Q.\end{equation}

In the classical field theory we need an affine version of the above procedure, since in the space of
infinitesimal configurations $\sJ^1E$ we have only an affine structure on fibres over $E$. Moreover, we have
to replace real-valued Lagrangians with Lagrangians taking values in the vector bundle $\Omega^m$.
\medskip

For simplicity, let us work first with an affine bundle $A\rightarrow N$ and a vector space $U$ like in the
previous section. The canonical evaluation between elements of $A$ and elements of $\Aff(A, U)$ is now a map
defined on a submanifold $A\times_N \Aff(A, U)\subset A\times \Aff(A, U)$ with values in $U$. The evaluation
in coordinates (chosen as in the previous section) reads
$$A\times_N \Aff(A, U)\ni (x^i, f^\alpha, \varphi_\alpha, r)\longmapsto (f^\alpha\varphi_\alpha+r)u\in U.$$
The canonical evaluation generates a Lagrangian submanifold of
$$\sT^\ast(A\times \Aff(A, U))\otimes_{(A\times \Aff(A, U))} U$$
with respect to the canonical $U$-valued symplectic form. There is an identification of the above cotangent
bundle with
$$\sT^\ast A\otimes_A U\times \sT^\ast\Aff(A, U)\otimes_{\Aff(A, U)} U.$$
 The canonical
symplectic form is identified with $\omega_A+\omega_{\Aff(A, U)}$. If we use the adapted coordinates $(x^i,
f^\alpha, \sigma_i, \psi_\alpha)$ in $\sT^\ast A\otimes_A U$ and appropriate coordinates $(x'{}^i,
\varphi_\alpha, r, \sigma'_i, h^\alpha,\rho)$ in $\sT^\ast\Aff(A, U)\otimes_{\Aff(A, U)} U$, we get the
generated Lagrangian submanifold given by the conditions:
$$x^i=x'{}^i,\qquad \sigma_i=-\sigma'_i, \qquad \psi_\alpha=\varphi_\alpha,\qquad f^\alpha=h^\alpha,\qquad \rho=1.$$
To get a symplectic relation out of that submanifold, we have to change signs in the fibre of
$\sT^\ast\Aff(A, U)\otimes_{\Aff(A, U)} U$. The graph of $\widetilde{ \mathcal{R}}_A$ is therefore given by
the conditions
$$x^i= x'{}^i,\qquad \sigma_i=\sigma'_i, \qquad \psi_\alpha=\varphi_\alpha,\qquad f^\alpha=-h^\alpha,\qquad \rho=-1.$$
We observe that the relation $\widetilde{ \mathcal{R}}_A$ is a map defined on the submanifold $K_{-id}$ with
values in $\sT^\ast A\otimes_A U$. In coordinates,
$$\widetilde{ \mathcal{R}}_A(x^i, \varphi_\alpha, r, \sigma_i, f^\alpha,-1)=(x^i, -f^a, \sigma_i, \varphi_a ).$$
Since $\widetilde{ \mathcal{R}}_A$ is constant on leaves of the characteristic foliation of $K_{-id}$, it
reduces to a symplectomorphism
$$\mathcal{R}_A: \sT^\ast A\otimes_A U\longrightarrow \sP \Aff(A, U).$$
For any function $L: A\rightarrow U$, from $\xd L(A)\subset \sT^\ast A\otimes_A U$ we obtain a Lagrangian
submanifold of $\sP \Aff(A, U)$. As a generating object we can choose a family of functions on $\Aff(A,U)$
parameterized by elements of $A$:
$$H: A\times_N\Aff(A,U)\longrightarrow U,\quad  H(a,\varphi)=L(a)-\varphi(a).$$
Note that $H$ satisfies condition (\ref{hm1}), i.e
$$H(a, \varphi+u)=H(a, \varphi)-u.$$
It generates, of course, a submanifold of $\sT^\ast\Aff(A, U)\otimes U$ which, after reduction, equals
$\mathcal{R}_A(\xd L(A)).$
If the family  reduces to a single function, the latter corresponds to a certain section $\Sigma_H$ of the
bundle $\theta:\Aff(A, U)\rightarrow V^\ast\otimes U$.
\smallskip

Applying again the above constructions to $A=\sJ^1_xE$ and recalling that we use the notation
$\Aff(\sJ^1_xE,\Omega^m_x)=\sJ^\dag_x E$, we get a diffeomorphism
$$\mathcal{R}_{\sJ^1 E}: \sV^\ast\sJ^1 E\otimes\Omega^m\longrightarrow \sP\sJ^\dag E$$
which, restricted to every fibre over $M$, is a symplectomorphism with respect to appropriate symplectic
$\Omega^m_x$ valued forms. We will denote with $\beta$ the composition
\begin{equation}\label{hm4}
  \beta=\alpha\circ\mathcal{R}_{\sJ^1 A},\qquad  \beta:\sJ^1\mathcal{P}\longrightarrow \sP\sJ^\dag E.
\end{equation}
In local coordinates, we get
\begin{equation}\label{be5}
\beta(q^i, y^a, p^j{}_b, y^c{}_k, p^l{}_{ds})=(q^i, y^a, p^j{}_b, \sum_k p^k{}_{ck}, y^d{}_{l}).
\end{equation}
The map $\beta$ constitutes the Hamiltonian side of the Tulczyjew triple for the classical field theory.
\begin{equation}
\xymatrix@C-5pt{
 & \sP\sJ^\dag E\ar[rd]^{\sP\varsigma}\ar[ldd]_{\sP\theta} & & & \sJ^1\mathcal{P}\ar[rd]^{\sJ^1\pi}\ar[ldd]_/-10pt/{\sj^1(\tau\circ\pi)}\ar[lll]_{\beta} &  \\
 & & \sJ^1 E\ar[ldd]_/-10pt/{\sj^1\zz} & & & \sJ^1 E\ar[ldd]_/-10pt/{\sj^1\zz}\ar[lll] \\
\mathcal{P}\ar[rd]^{\pi} & & & \mathcal{P}\ar[rd]^{\pi}\ar[lll] & & \\
 & E\ar[rrr]\ar[rrr] & & & E &
}
\end{equation}
{The above construction shows that, in the first-order field theory, a Hamiltonian is not a density valued
function on the phase space, but a section of certain affine bundle over phase space with one-dimensional
fibres. Differentials of such sections are elements of an affine analog of the cotangent bundle. From the
construction we obtain a family of Hamiltonian sections parameterized by elements of $\sJ^1 E$,
$$\Sigma_H: \sJ^1E\times_E\mathcal{P}\rightarrow \sJ^\dag E,$$
that corresponds to a family of density valued maps
$$H:\sJ^1 E\times_E \sJ^\dag E\rightarrow \Omega^m,\qquad H(\sj^1\sigma, \varphi)=L(\sj^1\sigma)-\varphi(\sj^1\sigma).$$
In some cases the above family reduces to a single generating section. It happens e.g. in Electrostatics (see
section \ref{secele}). In such cases we obtain $\mathcal{D}$ from the image of the differential of the
Hamiltonian section $H:\mathcal{P}\to\sJ^\dag E$ by means of the map $\beta$. More precisely,
\begin{equation}\label{dynham}\mathcal{D}=\beta^{-1}(\xd H(\mathcal{P})).\end{equation}
Like in the Lagrangian case, the process of generating phase dynamics from a Hamiltonian is very simple.}

\subsection{Structure of the phase space}

The phase space $\mathcal{P}$ is fibred over $E$. The fibration is a vector bundle. The space of vertical
vectors $\sV\mathcal{P}$ is therefore a double vector bundle fibrated over $\sV E$ and $\mathcal{P}$,
$$\xymatrix@C-5pt{
 & \sV \mathcal{P}\ar[dl]_{\nu_\mathcal{P}}\ar[dr]^{\sV\theta} & \\
\mathcal{P}\ar[dr]_\theta & & \sV E\ar[dl]^{\zn_E} \\
 & E &
}$$ We define a one-form $\vartheta_{\mathcal{P}}$ on $\mathcal{P}$ with values in $\Omega^{m-1}$ by the
formula
$$\vartheta_{\mathcal{P}}(\delta p)=\langle\nu_{\mathcal{P}}(\delta p), \sV\theta(\delta_p)\rangle.$$
In coordinates,  for $\delta p=(x^i, y^a, p^j{}_a, \delta y^c, \delta p^k{}_d)$, we get
$$\vartheta_{\mathcal{P}}(\delta p)=p^i{}_a\delta y^a\eta_i, \text{ i.e. }\vartheta_\sP= p^i{}_a\xd y^a\otimes\eta_i.$$
The form $\vartheta_\mathcal{P}$ is an analog of the canonical Liouville form on a cotangent bundle. Applying
$\xd_M$ to $\vartheta_\mathcal{P}$, we obtain a one-form on $\sJ^1 \mathcal{P}$ with values in $\Omega^m$
which in coordinates reads
$$\xd_M\vartheta_{\mathcal{P}}=p^i{}_{ai}\xd y^a\otimes\eta+ p^j{}_b\xd y^b{}_j\otimes \eta$$
and which is an analog of $\xd_\sT\vartheta_M$, the Liouville form on $\sT\sT^\ast M$.

Applying vertical differential $\xd$ to $\xd_M\vartheta_{\mathcal{P}}$, we obtain a two-form on $\sJ^1
\mathcal{P}$ which can be treated as a family of presymplectic forms with values in $\Omega^m$, parameterized
by points of $M$:
$$\omega_{\sJ^1\mathcal{P}}=\xd\xd_M\vartheta_{\mathcal{P}}=\xd p^i{}_{ai}\wedge\xd y^a\otimes\eta
+\xd p^j{}_b\wedge\xd y^b{}_j\otimes \eta.$$ It is easy to see in coordinates that
$$\beta^\ast\omega_{\sP\sJ^\dag E}=-\omega_{\sJ^1\mathcal{P}}.$$
There is an alternative construction of the map $\beta$ that uses the language of
differential forms on fiber bundles. The crucial role in the construction is played by the
two-form $\xd \vartheta_{\mathcal{P}}$ with values in $\Omega^{m-1}$ \cite{V2}.

\subsection{Electrostatics}\label{secele}

The Hamiltonian side of the Tulczyjew triple for Electrostatics is simplified, because fields are sections of
the trivial bundle $\zeta: A\times\R\rightarrow \R$. The bundle $\sJ^1\zeta$ is therefore a vector bundle:
$$\sJ^1\zeta: \sJ^1 E=A\times \R\times V^\ast\longrightarrow A\times \R.$$
The affine dual to this vector bundle is again the Cartesian product. Note that $\Omega^m\simeq
A\times\bigwedge^3 V^\ast$. We have the following identifications:
\begin{align*}
 \sJ^1E &\simeq A\times \R\times V^\ast, \\
 \sJ^\dag E & \simeq A\times \R\times \bigwedge^2{}_o V^\ast\times \bigwedge^3{}_o V^\ast,  \\
 \mathcal{P} & \simeq  A\times \R\times \bigwedge^2{}_o V^\ast.
\end{align*}
The bundle
$\theta: \sJ^\dag E\longrightarrow \mathcal{P}$
is trivial, i.e.
$$\theta: A\times \R\times \bigwedge^2{}_o V^\ast\times \bigwedge^3{}_o V^\ast\longrightarrow
A\times \R\times \bigwedge^2{}_o V^\ast.$$ There is a natural action of $\bigwedge^3{}_o V^\ast$ in the fibres
of the above bundle by addition:
\begin{eqnarray*}
&\left(\bigwedge^3{}_o V^\ast\right)\times\left( A\times \R\times \bigwedge^2{}_o V^\ast\times \bigwedge^3{}_o
V^\ast\right)
\longrightarrow  A\times \R\times \bigwedge^2{}_o V^\ast\times \bigwedge^3{}_o V^\ast, \\
&(u,(x,r,p,\lambda))\longmapsto(x,r,p,\lambda+u).
\end{eqnarray*}
The canonical evaluation between $(x,r,\mu)\in\sJ^1 E$ and $(x,r,p,\lambda)\in\sJ^\dag E$ with values in
$\bigwedge^3{}_oV^\ast$ reads
\begin{equation}
\langle (x,r,\mu), (x,r,p,\lambda)\rangle =\mu\wedge p+\lambda.
\end{equation}
The above evaluation generates a relation $\mathcal{R}$ between the cotangent bundle of the space of
infinitesimal configurations,
$$\sV^\ast \sJ^1 E\otimes\Omega^m\;\;\simeq\;\;
A\times \R\times V^\ast\times\bigwedge^3{}_o V^\ast\times\bigwedge^2{}_o V^\ast$$ and the cotangent bundle of
the affine dual,
$$\sV^\ast \sJ^\dag E\otimes\Omega^m\;\;\simeq\;\;
A\times \R\times \bigwedge^2{}_o V^\ast\times \bigwedge^3{}_o V^\ast \times V^\ast \times\bigwedge^3{}_o
V^\ast\times\R.$$
An element $(x,r,\mu,\sigma,\omega)\in\sV^\ast \sJ^1 E\otimes\Omega^m$
is in the relation $\mathcal{R}$ with an element
$(x,r,p,\lambda,a,b,c)\in \sV^\ast \sJ^\dag E\otimes\Omega^m,$
if and only if
$$p=\omega,\quad a=-\sigma,\quad b=\mu,\quad c=1.$$
There is an action of $\bigwedge^3 V^\ast$ in the cotangent bundle $\sV^\ast \sJ^\dag E\otimes\Omega^m$ lifted
from the action in $\sJ^\dag E$. The image of $\mathcal{R}$ (defined by $c=1$) is invariant with respect to
the lifted action and the quotient space
$$\sP\sJ^\dag E\quad\text{can be identified with}\quad A\times\R\times\bigwedge^2{}_oV^\ast\times \bigwedge^3{}_oV^\ast \times V^\ast. $$
The graph of $\mathcal{R}$ is also invariant with respect to the lifted action, therefore there exists the
quotient relation between $\sV^\ast \sJ^1 E\otimes\Omega^m$ and $\sP\sJ^\dag E$. This quotient relation is
actually a map,
\begin{eqnarray*}
& \tilde{\mathcal{R}}: \sV^\ast \sJ^1 E\otimes\Omega^m\longrightarrow \sP\sJ^\dag E, \\
& (x,r,\mu,\sigma,\omega)\longmapsto (x,r,\omega, -\sigma, \mu).
\end{eqnarray*}
Composing the map $\tilde{\mathcal{R}}$ with $\alpha$ from the Lagrangian side, we get the map $\beta$
\begin{equation*}\label{hmel1}
\beta: \sV^\ast E\otimes\Omega^m\times_E\sJ^1\mathcal{P}\longrightarrow\sP\sJ^\dag E,\\
\end{equation*}
i.e.
\begin{eqnarray*}\label{hmel2}
& A\times\R\times \bigwedge^3{}_{o}V^\ast\times \bigwedge^2{}_{o}V^\ast\times V^\ast\times
V^\ast\otimes\bigwedge^2{}_{o}V^\ast
\longrightarrow A\times\R\times\bigwedge^2_oV^\ast\times \bigwedge^3{}_oV^\ast \times V^\ast, \\
& (x,r,\rho,p,\mu,\nu)\longmapsto(x,r,p,\rho-Alt(\nu),\mu).
\end{eqnarray*}
Equation (\ref{elec3}) for a section of the phase bundle $\mathcal{P}\rightarrow A$ can be generated from a
Hamiltonian. The generating family
$$h(x,r,p, \lambda, \mu)=L(x,r,\mu)-\zm\wedge p-\lambda\mathfrak{g}$$
reduces to a section $H$ of the bundle $\theta$. Since the bundle is trivial, this section can be represented
by the map
$$H: \mathcal{P}\longrightarrow\bigwedge^3{}_oV^\ast,\qquad H(x,r,p)=\frac{1}{2}(\ast p)\wedge p.$$
where $\ast$ is the Hodge-star associated with the metric $g$. The Hamiltonian $H$ generates the subset
$$\xd H(\mathcal{P})\subset\sP\sJ^\dag E,\qquad \xd H(\mathcal{P})=\{(x,r,p,a,b): a=0, b=\ast p\}.$$
The inverse image of $\xd H(\mathcal{P})$ by $\beta$ is
\begin{equation}\label{hmel3}
\beta^{-1}(\xd H(\mathcal{P}))=\{(x,r,p,\rho,\mu,\nu):\; \rho=Alt(\nu), \mu=\ast p \}.
\end{equation}
Comparing (\ref{hmel3}) with (\ref{elect4}), we see that
$$D=\beta^{-1}(\xd H(\mathcal{P})).$$

\section{Conclusions}

We have constructed the Tulczyjew triple for the first-order field theory, starting from fundamental concepts of
calculus of variations. Our results can be summarized in the following diagram of affine and vector bundle morphisms
\begin{equation}\label{triple}
\xymatrix@!C=2pc{
 & \sP\sJ^\dag E\ar[dl]_{\sP\theta}\ar[ddr]^/-10pt/{\sP\varsigma} & & &
 \sJ^1\mathcal{P}\ar[lll]_\beta\ar[rrr]^\alpha\ar[dl]_{\sj^1(\tau\circ\pi)}\ar[ddr]^/-10pt/{\sJ^1\pi}& & &
 \sV^\ast\sJ^1E\otimes\Omega^m \ar[dl]_\xi\ar[ddr]^/-10pt/{\rho_{\sJ^1 E}}& \\
\mathcal{P}\ar[ddr]^/-10pt/{\pi} & & & \mathcal{P}\ar[rrr]\ar[lll]\ar[ddr]^/-10pt/{\pi} & & & \mathcal{P}\ar[ddr]^/-10pt/{\pi} & & \\
    & & \sJ^1 E\ar[dl]_{\sj^1\zz} & & & \sJ^1 E\ar[rrr]\ar[lll]\ar[dl]_{\sj^1\zz} & & & \sJ^1 E\ar[dl]_{\sj^1\zz} \\
 & E &  & & E\ar[rrr]\ar[lll] & & & E &
}
\end{equation}
{The left-hand side of the triple is Lagrangian, the right-hand side is Hamiltonian, and the phase dynamics
lives in the middle. The phase dynamics $\mathcal{D}$ being a subset of $\sJ^1\mathcal{P}$ is interpreted as a
condition for first jets of sections of the momentum bundle. It can be obtained from a Lagrangian as
$$\mathcal{D}=\alpha^{-1}(\xd L(\sJ^1 E))$$
or from a Hamiltonian as
$$\mathcal{D}=\beta^{-1}(\xd \Sigma_H(\mathcal{P})).$$
The Hamiltonian is a section of an affine bundle $\theta$ with fibres modeled on the vector spaces of volume
elements on $M$.}

 The Lagrangian side of the triple looks similarly in almost all papers devoted to
the Tulczyjew triple for field theory. However, we would like to emphasize that all the spaces and maps that
we use are constructed, not postulated, and have clear interpretation in the language of variational calculus.
It is interesting also that the phase space $\mathcal{P}$ is not dual to the space of infinitesimal
configurations $\sJ^1 E$.

To construct the Hamiltonian side of the triple we have used the notion of affine duality. Geometrical
language for the first-order theory is not the only place in classical mathematical physics where affine
structures are needed. Constructions that were used in section \ref{secham} are similar to those needed in
time-dependent mechanics \cite{GGU1} or in the intrinsic formulation of Newtonian mechanics \cite{GU}. We
expect also that the notion of affine duality will play an important role in higher-order theories. The
Hamiltonian formulation of the first-order field theory is based on the canonical isomorphism $\sP\sJ^\dag
E\simeq \sV^\ast\sJ^1E\otimes\Omega^m$ generated by the evaluation between $\sJ^1 E$ and its affine-dual
bundle.

In the Tulczyjew triple for mechanics, all three spaces, $\sT^\ast\sT^\ast Q$, $\sT\sT^\ast Q$, and $\sT^\ast\sT
Q$ are isomorphic. It is not the case in the field theory. The middle space is not isomorphic to $\sP\sJ^\dag E$
and $\sV^\ast\sJ^1E\otimes\Omega^m$. To have an isomorphism between all three spaces, we could replace the space
$\sJ^1\mathcal{P}$ with the quotient space with respect to a certain equivalence relation. We decided not to do
it, because passing to the quotient we would loose the obvious interpretation of the constitutive set as a first-order
differential equation. However, it is clear that only very special differential inclusions $\mathcal{D}\subset
\sJ^1 \mathcal{P}$ are generated by some Lagrangian or Hamiltonian.

The spaces $\sP\sJ^\dag E$ and $\sV^\ast\sJ^1E\otimes\Omega^m$ are equipped with canonical two
forms that, restricted to every fibre, are symplectic forms with values in $\Omega^m$. We
have also a canonical presymplectic forms on fibres of the bundle $\sJ^1\mathcal{P}\rightarrow M$. The phase
space $\mathcal{P}$ possesses a canonical one-form with values in $\Omega^{m-1}$ which is an analog of the
canonical Liouville form on $\sT^\ast M$.


\begin{thebibliography}{19}


\bibitem{B} S. Benenti, \emph{Hamiltonian optics and generating families},
Napoli Series on Physics and Astrophysics, {\bf 8} (2004), pp.218

\bibitem{CIL} F. Cantrijn, L. A. Ibort, M. De Leon, \emph{Hamiltonian structures on multisymplectic
manifolds}, Rend. Sem. Mat. Univ. Pol. Torino, {\bf 54} (1996), 225--236.

\bibitem{CCI1} J. F. Cari\~nena, M. Crampin, L. A. Ibort, \emph{On the multisymplectic formalism for frst order
theories}, Differential Geom. Appl., {\bf 1} (1991), 354--374.

\bibitem{CCI2} J. F. Cari\~nena, M. Crampin, L. A. Ibort, \emph{On the geometry of multisymplectic manifolds},
J. Austral. Math. Soc. A, {\bf 66} (1999), 303--330.

\bibitem{CGM} C. M. Campos, E. Guzm\'{a}n, J. C. Marrero, \emph{Classical field theories of first order and lagrangian
submanifolds of premultisymplectic manifolds}, arXiv:1110.4778.

\bibitem{EM} A. Echeverria-Enriquez, M. C. Mu\~noz-Lecanda, \emph{Geometry of multisymplectic Hamiltonian
first order theory}, J. Math. Phys., {\bf 41} (2000), 7402--7444.

\bibitem{FP1} M. Forger, C. Paufler, H. R\"{o}mer, \emph{A general construction of Poisson brackets on exact
multisymplectic manifolds}, Rep. Math. Phys., {\bf 51} (2003), 187--195.

\bibitem{FP2} M. Forger, C. Paufler, H. R\"{o}mer, \emph{Hamiltonian multivector fields and Poisson forms in
multisymplectic field theories,} J Math. Phys., {\bf 46} (2005), 112903--112932.

\bibitem{GM} G. Giachetta, L. Mangiarotti, \emph{Constrained Hamiltonian Systems and Gauge Theories},
{ Int. J. Theor. Phys.}, {\bf 34} (1995), 2353--2371.

\bibitem{GIM1} M. J. Gotay, J. Isenberg, J. E. Marsden, \emph{Momentum maps and classical relativistic fields, Part I: Covariant
field theory}, { arXiv:physics/9801019v2} (2004).

\bibitem{GIM2} M. J. Gotay, J. Isenberg, J.E. Marsden, \emph{Momentum maps and classical relativistic fields, Part II:
Canonical analysis of field theories.} {arXiv:math-ph/0411032v1} (2004).

\bibitem{G1} M. J. Gotay, \emph{A multisymplecitc framework for classical field theory and the calculus of variations I:
Covariant Hamiltonian Formalism}, { in Francaviglia, M., editor, Mechanics, Analysis, and Geometry: 200 Years
After Lagrange}, North Holland, Amsterdam, (1991) 203–-235.

\bibitem{G2} M. J. Gotay, \emph{A multisymplectic framework for classical field theory and the calculus of vari-
ations II. Space + time decomposition}, Differential Geom. Appl., {\bf 1} (1991), 375--390.

\bibitem{G} K. Grabowska, \emph{Lagrangian and Hamiltonian formalism in field theory: a simple model}, {J. Geom.
Mech.}, {\bf 2} (2010), 375--395.

\bibitem{GU} K. Grabowska, P. Urba\'nski: \emph{AV-differential geometry and Newtonian Mechanics}, {Rep. Math. Phys.}
{\bf 58} (2006) 21--40.

\bibitem{GGU1} K. Grabowska, J. Grabowski, P. Urba\'nski,
\emph{AV-differential geometry: Poisson and Jacobi structures},  J. Geom. Phys., {\bf 52} (2004), 398--446. 

\bibitem{GGU0} K. Grabowska, J. Grabowski, P. Urba\'nski,
\emph{AV-differential geometry: Euler-Lagrange equations},  J. Geom. Phys., {\bf 57} (2007), 1984–-1998 

\bibitem{GGU2} K. Grabowska, J. Grabowski, P. Urba\'nski, \emph{Geometrical Mechanics on algebroids},
Int. J. Geom. Meth. Mod. Phys., {\bf 3} (2006), 559--575. 

\bibitem{GR} J.~Grabowski, M.~Rotkiewicz, \emph{Higher vector bundles and multi-graded symplectic
manifolds},  J. Geom. Phys., {\bf 59} (2009), 1285-1305.

\bibitem{GRU} J.~Grabowski, M.~Rotkiewicz, P.~Urba\'nski,
\emph{Double affine bundles}, {J. Geom. Phys.} {\bf 60} (2010), 581-598.

\bibitem{HK} F. Helein, J. Kouneiher, \emph{Covariant Hamiltonian formalism for the calculus of variations
with several variables: Lepage-Dedecker vs. De Donder-Weyl}, Adv. Theor. Math. Phys., {\bf 8} (2004),
565--601.

\bibitem{Kij} J. Kijowski, \emph{Elasticit\`a finita e relativistica: introduzione
ai metodi geometrici della teoria dei campi}, Pitagora Editrice (Bologna) (1991).

\bibitem{Kij2} J. Kijowski, W. M. Tulczyjew, \emph{A Symplectic Framework for Field Theories}, Lecture Notes in Physics,
{\bf 107}, Springer-Verlag, Berlin-New York, (1979).

\bibitem{KU} K.~Konieczna, P.~Urba\'nski, \emph{Double vector
bundles and duality}, Arch. Math. (Brno), {\bf 35} (1999), 59--95.

\bibitem{Kr1} O. Krupkov\'{a}, \emph{Hamiltonian field theory}, J. Geom. Phys., {\bf 43} (2002), 93--132.

\bibitem{Kr2} O. Krupkov\'{a}, P. Voln\'{y}, \emph{Euler-Lagrange and Hamilton equations for non-holonomic systems in field theory},
J. Phys. A: Math. Gen. {\bf 38} (2005), 8715--874.

\bibitem{Kr3} O. Krupkov\'{a}, P. Voln\'{y}, \emph{Differential equations with constraints in jet bundles: Lagrangian and Hamiltonian systems},
Lobachevskii Journal of Mathematics {\bf 23} (2006), 95--150.

\bibitem{LMS} M. de Le\'on, D. Mart\'in de Diego, A. Santamar\'ia-Merino,
\emph{Tulczyjew's triples and lagrangian submanifolds in classical field theories}, in ``Applied Differential
Geometry and Mechanics," Editors W. Sarlet and F. Cantrijn, Univ. of Gent, Gent, Academia Press (2003),
21–-47.

\bibitem{MTU} G. Marmo, W.M. Tulczyjew, P. Urba\'nski, \emph{Dynamics of Autonomous Systems with External Forces},
\emph{Acta Physica Polonica B}, {\bf 33}, No. 5, (2002), p. 1181.

\bibitem{Pr1} J.~Pradines, \emph{Fibr\'es vectoriels doubles et calcul
des jets non holonomes} (French), Notes polycopi\'ees, Amiens, (1974).

\bibitem{RRR} A. M. Rey, N. Roman-Roy, M. Salgado, S. Vilari\~{n}o, \emph{On the $k$-symplectic, $k$-cosymplectic and multisymplectic formalism of classical field theories}, J. Geom. Mech. {\bf 3} (2011), 113-–137.

\bibitem{ST} J. Sniatycki, W.M. Tulczyjew, \emph{Generating Forms of Lagrangian Submanifolds},
Indiana University Mathematics Journal, {\bf 22}, (1972) pp. 267--275.

\bibitem{Tu8} M.R. Menzio, W.M. Tulczyjew, \emph{Infinitesimal symplectic
relations and generalized Hamiltonian dynamics}, Ann. Inst. H. Poincare, {\bf 28}, (1978), pp. 349--367.

\bibitem{Tu7} W.M. Tulczyjew, \emph{Relations symplectiques et les equations d'Hamilton-Jacobi relativistes},
C.R. Acad. Sc. Paris, {\bf 281}, (1975), pp. 545--548.

\bibitem{Tu9} W.M. Tulczyjew, \emph{Les sous-varietes Lagrangiennes et la Dynamique Hamiltonienne},
C.R. Acad. Sc. Paris, {\bf 283}, (1976), pp. 15--18.

\bibitem{Tu10} W.M. Tulczyjew, \emph{Les sous-varietes Lagrangiennes et la Dynamique Lagrangienne},
C.R. Acad. Sc. Paris, {\bf 283}, (1976), pp. 675--678.

\bibitem{Tu5} W.M. Tulczyjew, \emph{A symplectic framework for linear field theories},
Annali di Matematica pura ed applicata, {\bf 130}, (1982), pp. 177--195.

\bibitem{Tu6} W.M. Tulczyjew, \emph{The Legendre Transformation},
Ann. Inst. H. Poincare, {\bf 27}, (1977), pp. 101--114.

\bibitem{Tu1} W. M.~Tulczyjew, \emph{Hamiltonian
systems, Lagrangian systems, and the Legendre transformation}, Symposia Mathematica, Vol. XIV (Convegno di
Geometria Simplettica e Fisica Matematica, INDAM, Rome, 1973),
pp. 247-–258. Academic Press, London, (1974). 

\bibitem{TU3}  W. M. Tulczyjew, \emph{Geometric Formulation of Physical Theories}, Bibliopolis, (1989).

\bibitem{Tu2}  W. M. Tulczyjew, \emph{The Euler-Lagrange resolution}, in
Differential Geometrical Methods in Mathematical Physics, Lecture Notes in Mathematics, {\bf 836}, (1980),
22--48.

\bibitem{TU}  W. M.~Tulczyjew, P. Urba\'nski,  \emph{A slow
and careful Legendre transformation for singular Lagrangians}, The Infeld Centennial Meeting (Warsaw,
1998), Acta Phys. Polon. B,  {\bf 30} (1999), 2909--2978. 

\bibitem{TU2} W. M. Tulczyjew, P. Urba\'nski, \emph{Liouville structures},
Universitatis Iagellonicae Acta Mathematica {\bf 47} (2009), pp. 187--226.

\bibitem{VCL} J. Vankershaver, F. Cantrijn, M. De Leon and M. De Diego, \emph{Geometric aspects of nonholonomic field theories}, Rep. Math. Phys., {\bf 46} (2005), 387--411.

\bibitem{V} L. Vitagliano, \emph{The Hamilton-Jacobi Formalism for Higher Order Field Theories},
Int. J. Geom. Meth. Mod. Phys, {\bf 07} (2010) pp.1413--1436.

\bibitem{V2} L. Vitagliano, \emph{private communication}.

\end{thebibliography}
\end{document}